\documentstyle[amssymb,amsfonts]{amsart}

\newcommand{\E}{{\cal E}}

\def\R{{\mbox{\bf R}}}

\def\Om{\Omega}

\def\eps{\varepsilon}
\def\pp{{p^\prime}}
\def\np{{n^\prime}}
\def\qp{{q^\prime}}
\def\kp{{k^\prime}}
\def\x{{\underline{x}}}
\newcommand{\qtil}{{\tilde{q}}}
\renewcommand{\mp}{{m^\prime}}
\newcommand{\ktil}{{\tilde{k}}}
\newcommand{\ptil}{{\tilde{p}}}
\def\emph#1{{\it #1}}
\def\textbf#1{{\bf #1}}
\newenvironment{proof}{\noindent {\bf Proof} }{\endprf\par}
\def \endprf{\hfill  {\vrule height6pt width6pt depth0pt}\medskip}

\theoremstyle{plain}
  \newtheorem{theorem}[subsection]{Theorem}
  \newtheorem{proposition}[subsection]{Proposition}
  \newtheorem{lemma}[subsection]{Lemma}
  \newtheorem{corollary}[subsection]{Corollary}
  \newtheorem{conjecture}[subsection]{Conjecture}

\theoremstyle{remark}

\theoremstyle{definition}
  \newtheorem{definition}[subsection]{Definition}

\include{psfig}

\begin{document}

\title[Bilinear restriction and Kakeya]{A bilinear approach
to the restriction and Kakeya conjectures}

\author{Terence Tao}
\address{Department of Mathematics, UCLA, Los Angeles, CA 90024}
\email{tao@@math.ucla.edu}

\author{Ana Vargas}
\address{Departamento de Matem\'aticas, Universidad Aut\'onoma
de Madrid, 28049 Madrid (Spain).}
\email{ana.vargas@@uam.es}

\author{Luis Vega}
\address{Departamento de Matem\'aticas, Universidad del
Pa\'is Vasco, Apartado 644, 48080, Bilbao (Spain).}
\email{mtpvegol@@lg.ehu.es}

\subjclass{42B10, 42B25}
\keywords{Restriction conjecture, bilinear estimates, Kakeya conjecture}
\begin{abstract}
Bilinear restriction estimates have been appeared in work of Bourgain,
Klainerman, and Machedon.  In this paper we develop the theory of these 
estimates (together with
the analogues for Kakeya estimates).
As a consequence we improve the $(L^p,L^p)$ spherical restriction theorem
of Wolff \cite{wolff:kakeya} from $p > 42/11$ to $p > 34/9$, and also
obtain a sharp $(L^p,L^q)$ spherical restriction theorem 
for $q> 4 - \frac{5}{27}$.
\end{abstract}

\maketitle

\tableofcontents

\section{Introduction}

The purpose of this paper is to investigate bilinear variants of the
restriction and Kakeya conjectures, to relate them to the standard
formulations of these conjectures, and to give applications of this
bilinear approach to existing conjectures.  The methods used are
based on several observations and results of Bourgain (see
\cite{borg:kakeya}-\cite{borg:stein}),
together with some refinements by Moyua, Vargas, and Vega 
\cite{vargas:restrict, vargas:2}.

This paper is organized as follows.  In the first section we discuss
bilinear restriction estimates, and show how one can pass back and forth
between these estimates and the standard restriction estimates.  We
also generalize the $12/7$ bilinear restriction estimate of 
\cite{vargas:2} to higher dimensions.

In the second section we give analogues of the above results
for the Kakeya operator.  In particular we give a bilinear improvement
to Wolff's Kakeya theorem in arbitrary dimension.

In the third section we give applications of these bilinear estimates in three
dimensions.  For 
example, we are able to improve the $42/11$ exponent in Wolff's restriction
theorem to $34/9$.  We are also able to prove a sharp $(L^p,L^q)$ restriction
theorem which improves on the classical $(L^2, L^4)$
Tomas-Stein theorem, and also give some concrete progress on a
bilinear restriction conjecture of Klainerman and Machedon.  We also give a 
non-bilinear approach 
to these estimates, which gives weaker results but is more direct and 
probably has a wider
range of application. 

Finally, we collect some standard harmonic analysis estimates in an
Appendix for easy reference.

This work was conducted at MSRI (NSF grant 9701955).  The authors wish 
to thank Tony Carbery, Adela Moyua, and Wilhelm Schlag for many helpful 
discussions.  The second author was partially supported by
the Spanish DGICYT (grant number PB94-149)
and the European Comission via the TMR network (Harmonic Analysis).

\section{Bilinear restriction estimates}\label{rest-sec}

Fix\footnote{All constants in this section are assumed to depend only on
$n$ and $A$.} $n \geq 2$ and $A > 0$,
and let $Q$ be the cube $[-1,1]^{n-1}$ in $\R^{n-1}$.
Let $\Phi: Q \to \R$ be a phase function satisfying the following
conditions:
\begin{itemize}
\item $\|\partial^\alpha \Phi\|_\infty \leq A$ for all
$0 \leq |\alpha| \leq N$, where $N$ is a large constant.
\item $\Phi(0) = \nabla \Phi(0) = 0$.
\item For all $x \in Q$,
the eigenvalues of the Hessian $\Phi_{x_i x_j}(x)$ all lie in
$[1-\epsilon_0, 1+\epsilon_0]$, where $\epsilon_0 > 0$ is a small constant.
\end{itemize}
We will call such a phase \emph{elliptic}.
The model example of an elliptic phase function is of course
the quadratic phase
$\Phi(x) = \frac{1}{2}|x|^2$, but any smooth
compact convex surface can be decomposed into finitely many graphs
whose graphing function (after an affine transformation) obeys the
above properties.  In particular, the unit sphere can be decomposed
in this manner.

We will consider linear and bilinear bounds for the
operator $\Re^*: L^1(Q) \to L^\infty(\R^n)$ defined by
$$ \Re^* f(\underline{x}, x_n) = \int_Q e^{-2\pi i (\underline{x} \cdot y
+ x_n \Phi(y))} f(y)\ dy.$$
This operator can be thought of as an adjoint restriction operator
associated to the surface $\{ (y, \Phi(y)): y \in Q\}$.
For $0 < p, q \leq \infty$, we use $R^*(p \to q)$ to denote
the estimate
$$ \| \Re^* f \|_q \lesssim \|f\|_p$$
for all test functions $f$, with the constant depending only on
$n$ and $A$.  Similarly, we use $R^*(p_1 \times p_2 \to q)$
to denote the estimate
$$ \| \Re^* f \Re^* g\|_q \lesssim \|f\|_{p_2} \|g\|_{p_1},$$
for all test functions $f$, $g$ supported on $Q_1$, $Q_2$
respectively, where $Q_1$, $Q_2$ are any sub-cubes of $Q$ whose
size and separation are comparable to 1.  (We will call such
cubes $O(1)$-separated in the sequel).

Estimates of the form $R^*(p \to q)$ are adjoint restriction estimates
and have attracted wide interest.  The (sharp) restriction conjecture
states that 
\begin{conjecture}\label{rest} $R^*(p \to q)$ holds 
whenever $q > \frac{2n}{n-1}$
and $p^\prime \leq \frac{n-1}{n+1} q$.
\end{conjecture}
These conditions are
well known to be best possible (see e.g. \cite{tomas:restrict}).  
This conjecture has been verified for $n=2$ \cite{carl:disc}, but
remains open in higher dimensions.  The main difficulty lies in making
the $q$ exponent as low as possible; the estimate is trivial for $q=\infty$,
H\"older's inequality can be used to raise $p$, and in certain cases
factorization theory can be used to lower $p$.  When $(p,q)$ lie
on the sharp line
\begin{equation}\label{sharp}
p^\prime = \frac{n-1}{n+1} q
\end{equation}
we abbreviate the estimate $R^*(p \to q)$ 
to $R^*_s(q)$.

We summarize\footnote{Some of the earlier results were not
stated for arbitrary elliptic phase functions.}
the known results in $n=3$ in Table 1.  The classical
theorem of Tomas and Stein states that $R^*(2 \to \frac{2(n+1)}{n-1})
=R^*_s(\frac{2(n+1)}{n-1})$
for any $n \geq 2$.  Later improvements have been made on this
result \cite{borg:kakeya},\cite{borg:stein},\cite{wolff:kakeya}; in
particular, Moyua, Vargas, and Vega \cite{vargas:restrict, 
vargas:2} have recently observed
that one has the estimate $R^*(\frac{7}{3}+\eps \to 
\frac{42}{11}+\eps)$ in three dimensions.  However, none of these improvements
to the Tomas-Stein theorem lies on the sharp line
$p^\prime = \frac{n-1}{n+1} q$.  As one of the applications
of this paper we will prove a new restriction theorem on this
sharp line.

\begin{table}
\begin{tabular}{lll} 
1.& $R^*(1\to\infty)=R^*_s(\infty)$& Riemann-Lebesgue \\
2.& $R^*(2 \to 6)$& Stein, 1967 \cite{feff:thesis}\\
3.& $R^*(2 \to 4+\eps)$ & Tomas, 1975 \cite{tomas:restrict}\\
4.& $R^*(2 \to 4) = R^*_s(4)$ & Stein, 1975 (For $n=3$: Sj\"olin, 1972) \\
5.& $R^*(4-\frac{2}{15}+\eps \to 4-\frac{2}{15}+\eps)$
& Bourgain, 1991 \cite{borg:kakeya}\\
6.& $R^*(4-\frac{2}{11}+\eps \to 4-\frac{2}{11}+\eps)$
& Wolff, 1995 \cite{wolff:kakeya} \\
7.& $R^*(\frac{7}{3}+\eps \to 4-\frac{2}{11}+\eps)$ & Moyua, Vargas,
Vega, 1995 \cite{vargas:restrict, vargas:2} \\
8.&$R^*(\frac{170}{77}+\eps \to 4-\frac{2}{9}+\eps)$& Theorem \ref{main}\\
9.&$R^*_s(4-\frac{5}{27}+\eps)$& Theorem \ref{main}\\
?.&$R^*_s(3+\eps)?$ & (critical value)\\
\end{tabular}

\caption{Known restriction theorems for $n=3$.  $\eps$ denotes an arbitrary
positive number.}
\end{table}

Our improvements will be based on
 the bilinear restriction estimates 
$R^*(p \times p \to q)$ defined earlier, which we will now discuss.  
These estimates have appeared implicitly in many works (e.g. \cite{borg:cone},
\cite{vargas:2}), and are closely 
related to null form estimates for the wave equation (see 
\cite{kl-mac:cpde-null}-\cite{kl-mac:duke-null}; related ideas also
appear in \cite{beals:xsb}), but do not appear to have been 
explicitly studied until very recently.

When $(p, 2q)$ lie in the range predicted by Conjecture \ref{rest}
then $R^*(p,2q)$ and $R^*(p \times p \to q)$ are almost equivalent.  Indeed,
in Section \ref{iff} we will prove

\begin{theorem}\label{imply} Let $n \geq 2$ and $1 < p,q < \infty$ be such that
$2q > \frac{2n}{n-1}$ and $p^\prime \leq \frac{n-1}{n+1} 2q$.

Then $R^*(p, 2q)$ implies $R^*(p \times p \to q)$.  Furthermore,
if $R^*(\ptil \times \ptil \to \qtil)$
holds for all $(\frac{1}{\ptil},\frac{1}{\qtil})$ in a 
neighbourhood of $(\frac{1}{p},\frac{1}{q})$,
then $R^*(p, 2q)$ holds.
\end{theorem}

However, the bilinear estimate $R^*(p \times p \to q)$ can hold for
exponents which are not covered by the above theorem.  For instance, when 
$n=2$ an easy computation using Plancherel's theorem and a 
change of variables shows that 
$R^*(2 \times 2 \to 2)$
holds, even though the Knapp example shows that $R^*(2,4)$ fails
completely. 
Thus one expects the range of exponents for the bilinear restriction
estimate to be larger than that of Conjecture \ref{rest}.
For $n=3$ the first results in this direction were by 
Bourgain \cite{borg:cone} (although the theorem
$R^*(\frac{16}{9} \times \frac{16}{9} \to 2)$ implicitly
appeared in \cite{borg:16-9});
more recently, Moyua, Vargas, and Vega \cite{vargas:2} showed
that
\begin{equation}\label{12-7}
R^*(\frac{12}{7} \times \frac{12}{7} \to 2)
\end{equation}
for $n=3$.  We modestly generalize this result to higher dimensions as

\begin{theorem}\label{modest}  Suppose that $n \geq 2$.  Then
$$R^*(p \times p \to 2)$$
holds if and only if $p \geq \frac{4n}{3n-2}$.
\end{theorem}

Recently\footnote{Workshop in Harmonic Analysis and PDE, MSRI, July 1997.} 
Klainerman and Machedon conjectured that 
\begin{equation}\label{kl}
R^*(2 \times 2 \to \frac{n+2}{n})
\end{equation}
for all $n \geq 2$.  By interpolating \eqref{kl} with what is implied by
Conjecture \ref{rest}, one is led to the following

\begin{conjecture}\label{conj}
If $n \geq 2$, then $R^*(p \times p \to q)$ holds whenever
\begin{align}
q &\geq \frac{n}{n-1},\label{c0}\\
\frac{n+2}{2q} + \frac{n}{p} &\leq n\label{c1}\\
\frac{n+2}{2q} + \frac{n-2}{p} &\leq n-1.\label{c2}
\end{align}
\end{conjecture}

By Theorem \ref{modest}
and interpolation the conjecture is verified for $q \geq 2$ (and thus
for $n = 2$).
The exponents in the above
conjecture are best possible; we will sketch the proof of this
statement in Section \ref{nec}.  From Theorem \ref{imply} we see that
Conjecture \ref{conj} implies Conjecture \ref{rest}.

We depict the conjectured ranges for the estimates
$R^*(p \to 2q)$ and $R^*(p \times p \to q)$ in Figure 1.  The
restriction conjecture states that $R^*(p \to 2q)$ 
holds for all $(p,q)$ in the
trapezoidal region bounded by $1$, $c$, $d$, and $0$, except
for the upper line between $c$ and $d$ inclusive; by the above
Theorem, this is almost equivalent to
$R^*(p\times p \to q)$ holding in this region.  
Klainerman's conjecture asserts that 
$R^*(p \times p \to q)$ holds at the endpoint $b$.  The combined
Conjecture \ref{conj} states that $R^*(p \times p \to q)$
holds in the pentagonal region bounded by $1$, $b$, $c$, $d$, and
$0$, including the upper line mentioned previously;
this region is best possible.

By Theorem \ref{imply} the standard restriction estimate
$R^*(p\to 2q)$ and the bilinear estimate $R^*(p \times p \to q)$
are essentially equivalent in the line between $c$ and $1$.
The points
$1-7$ correspond to the standard restriction results,
while the point
$a$ corresponds to the bilinear restriction theorem \eqref{12-7}.

From Theorem \ref{imply} and bilinear interpolation 
it is possible
to obtain new linear and bilinear restriction theorems; for instance,
by interpolating between the bilinear form\footnote{i.e. $R^*(\frac73\times\frac73,\frac{21}{11})$.} of the result in 
\cite{vargas:restrict} and \eqref{12-7} and using Theorem \ref{imply},
one may obtain the sharp restriction theorem $R^*_s(q)$ for all
$q > 4 - \frac{2}{17}$.  We will improve on these results in
in Section \ref{apply}.

\begin{figure}[htbp] \centering
\ \psfig{figure=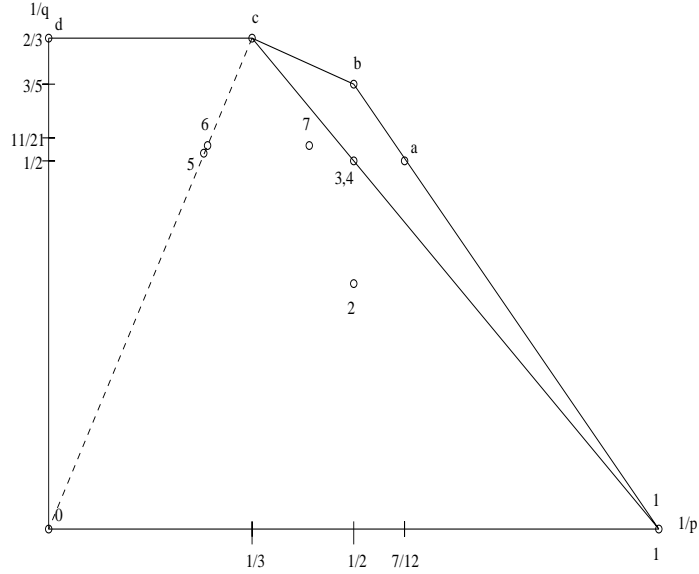,height=3in,width=3.6in}
\caption{Prior status of $R^*(p \times p \to q)$ and $R^*(p \to 2q)$ for $n=3$.  }
        \label{fig:bilinear}
	\end{figure}

\subsection{Proof of Theorem \ref{imply}}\label{iff}
\

The first implication is a trivial consequence of H\"older's inequality,
so we concentrate on the latter.  In view of the known results
for $n=2$ we may take $n \geq 3$.  From the Tomas-Stein theorem
(see e.g. \cite{stein:large}) and the necessity of \eqref{c0}
it suffices to consider the case $\frac{2n}{n-1} < 2q < 
\frac{2(n+1)}{n-1}$.  In particular we may assume that $1 < q < 2$.

The bilinear hypothesis $R^*(\ptil \times \ptil \to \qtil)$ allows us
to control $\Re^* f \Re^* g$ if $f$ and $g$ have $O(1)$-separated
supports.  By a parabolic rescaling argument this will imply a similar
estimate when $f$ and $g$ have $O(2^{-j})$-separated supports for any
$j > 0$.  Piecing these estimates together one may obtain an estimate
on $\Re^* f \Re^* g$ for arbitrary $f$, $g$, from which the conclusion
$R^*(p \to q)$ will follow.

We now turn to the details.  
Assume that the hypotheses of Theorem \ref{imply} hold.  We have to show that
$$ \|\Re^* f\|_{2q} \lesssim \|f\|_{p}.$$
By Marcinkeiwicz interpolation it suffices to show the restricted estimate
$$ \|\Re^* \chi_\Omega\|_{2q} \lesssim |\Omega|^{1/p}$$
for a slightly better value of $(p, 2q)$, where $\Omega$ is some
arbitrary subset of $Q$.

Let $j_0$ be the positive integer such that $|\Omega| \sim 2^{-j_0(n-1)}$.
Then by squaring the above estimate, we reduce
ourselves to 
\begin{equation}\label{target-1}
 2^{\frac{2(n-1)}{p}j_0} 
 \| \Re^* \chi_\Omega \Re^* \chi_\Omega \|_{q} 
 \lesssim 1.
\end{equation}
The next step is a Whitney decomposition.
For each $j > 0$, we dyadically decompose $Q$ into 
$\sim 2^{(n-1)j}$ dyadic subcubes $\tau^j_k$ of sidelength $2^{-j}$ in
the usual manner.  If $\tau^j_k$, $\tau^j_{\kp}$ are
two cubes with the same sidelength which are not adjacent but have
adjacent parents, we say that these cubes are \emph{close} and write
$\tau^j_k \sim \tau^j_{\kp}$.  For almost every $x, y \in Q$ there exists 
a unique pair of close cubes $\tau^j_k$, $\tau^j_{\kp}$ containing
$x$ and $y$ respectively.  Thus we have
$$ \Re^* \chi_\Omega \Re^* \chi_\Omega = \sum_j \sum_{k,\kp: \tau^j_k \sim \tau^j_{\kp}}
\Re^* \chi_{\Omega \cap \tau^j_k} \Re^* \chi_{\Omega \cap \tau^j_{\kp}}
$$
Thus to prove
\eqref{target-1} it suffices to show that
\begin{equation}\label{target-2}
2^{\frac{2(n-1)}{p}j_0}
\| \sum_{k, \kp: \tau^j_k \sim \tau^j_{\kp}}
\Re^* \chi_{\Omega \cap \tau^j_k} \Re^* \chi_{\Omega \cap \tau^j_{\kp}}\|_q
\lesssim 2^{-\eps |j - j_0|}
\end{equation}
for all $j > 0$ and some $\eps > 0$, since \eqref{target-1} 
follows from the triangle inequality.  
Informally, the above estimate asserts that the most significant separation 
scale is of the order of $2^{-j_0} = |\Omega|^{1/(n-1)}$;
this is already 
evident from the Knapp example.

Our next reduction will be to exploit some quasi-orthogonality between
the functions $\Re^* \chi_{\Omega \cap \tau^j_k} \Re^* \chi_{\Omega \cap \tau^j_{\kp}}$.  From the definition of $\Re^*$ we see that the Fourier transform
$\Re^* \chi_{\Omega \cap \tau^j_k}$ is supported on the infinite tube
$\tau^j_k \times \R$.  Thus,
the Fourier transform of
$\Re^* \chi_{\Omega \cap \tau^j_k} \Re^* \chi_{\Omega \cap \tau^j_{\kp}}$
is supported in the tube
$$ T_{j,k} = \tilde \tau^j_k \times \R,$$
where $\tilde \tau^j_k$ is a cube of sidelength $C 2^{-j}$ whose center
is twice that of $\tau^j_k$.  From Lemma \ref{quasi-lemma} in the Appendix
and the assumption $q < 2$, we thus have
$$
\|\Re^* \chi_{\Omega \cap \tau^j_k} \Re^* \chi_{\Omega \cap \tau^j_{\kp}}\|_q
\lesssim 
(\sum_k \|\Re^* \chi_{\Omega \cap \tau^j_k} \Re^* \chi_{\Omega \cap \tau^j_{\kp}}\|_q^q)^{1/q}.
$$
Thus \eqref{target-2} will be proven if we can show that
\begin{equation}\label{target-3}
2^{\frac{2(n-1)q}{p}j_0}
\sum_k \| \sum_{\kp: \tau^j_k \sim \tau^j_{\kp}}
\Re^* \chi_{\Omega \cap \tau^j_k} \Re^* \chi_{\Omega \cap \tau^j_{\kp}}\|_q^q
\lesssim 2^{-\eps q |j - j_0|}
\end{equation}

This will follow from the following estimate.

\begin{proposition}\label{rescale} 
For all $\ptil$ in a neighbourhood of $p$, we have
\begin{equation}\label{target-4}
\| \Re^* \chi_{\Omega \cap \tau^j_k} \Re^* \chi_{\Omega \cap \tau^j_{\kp}}\|_q
\lesssim 
2^{-\frac{2(n-1)}{\ptil^\prime}j} 2^{\frac{n+1}{q}j}
|\Omega\cap \tau^j_k|^{1/\ptil} |\Omega\cap \tau^j_\kp|^{1/\ptil}. 
\end{equation}
\end{proposition}

\begin{proof}
This will be accomplished
by a parabolic
rescaling argument.  By translating $\Phi$ and subtracting a 
harmless affine factor\footnote{Strictly speaking, one may need
to increase $A$ by a constant factor to do this; we will gloss over
this technicality.} we may assume that $\tau^j_k$ is centered
at the origin.  We now observe that since $\Phi$ is an elliptic
phase, the function 
$$ \tilde \Phi(x) = 2^{2j} \Phi(2^{-j} x)$$
is also elliptic.  
Since $R^*(\ptil \times \ptil \to q)$
holds for all $\ptil$ in a neighbourhood of $p$ by assumption, we have
$$ \| \tilde{\Re}^* f \tilde{\Re}^* g \|_q \lesssim \|f\|_\ptil \|g\|_\ptil$$
whenever $f$ and $g$ are supported on disjoint $O(1)$-separated cubes, where
$\tilde \Re^*$ is the adjoint restriction operator corresponding to $\tilde
\Phi$.  Applying a parabolic scaling $(\underline{x}, x_n)
\to (2^j \underline{x}, 2^{2j} x_n)$ to this estimate one obtains
$$ \| \Re^* f \Re^* g \|_q \lesssim 
2^{-\frac{2(n-1)}{\ptil^\prime}j} 2^{\frac{n+1}{q}j}
\|f\|_\ptil \|g\|_\ptil$$
whenever $f$ and $g$ are supported on $\tau^j_k$ and $\tau^j_\ktil$ 
respectively, and \eqref{target-4} follows.
\end{proof}

It remains to obtain \eqref{target-3} from the proposition.  Let
$\ptil < p$ be such that \eqref{target-4} holds.
If we apply \eqref{target-4}
and the triangle inequality,
we see that \eqref{target-3} reduces to
$$ 2^{\frac{2(n-1)q}{p}j_0} 
\sum_k (\sum_{\kp: \tau^j_k \sim \tau^j_{\kp}}
2^{-\frac{2(n-1)}{\ptil^\prime}j} 2^{\frac{n+1}{q}j}
|\Omega\cap \tau^j_k|^{1/\ptil} |\Omega\cap \tau^j_\kp|^{1/\ptil})^q
\lesssim 2^{-\eps q |j-j_0|}.$$
By polarization and the fact that for each $k$ there are only finitely
many cubes $\tau^j_{\kp}$ close to $\tau^j_k$, this is in turn reduces
to
\begin{equation}\label{major}
 2^{\frac{2(n-1)q}{p}j_0} 
 2^{-\frac{2(n-1)q}{\ptil^\prime}j} 2^{(n+1)j}
 \sum_k
|\Omega\cap \tau^j_k|^{\frac{2q}{\ptil}} \lesssim 2^{-\eps q |j-j_0|}.
\end{equation}
We divide into two cases: $\ptil \leq 2q$ and $\ptil > 2q$.  If
$\ptil \leq 2q$ then we use \eqref{xr-big} from Lemma \ref{xr-est}
in the Appendix with $\alpha=1$ to obtain
$$ \sum_k
|\Omega\cap \tau^j_k|^{\frac{2q}{\ptil}} \lesssim
2^{-(n-1)j_0} 2^{-(n-1)\max(j,j_0)(\frac{2q}{\ptil}-1)}.$$
Thus \eqref{major} reduces to
\begin{equation}\label{exp}
\begin{split}
\frac{2(n-1)q}{p}j_0  
-\frac{2(n-1)q}{\ptil^\prime}j& + (n+1)j
-(n-1)j_0\\
&-(n-1) \max(j,j_0) (\frac{2q}{\ptil}-1)
\leq -\eps q |j - j_0|.
\end{split}
\end{equation}
By convexity it suffices to verify this inequality for the values
$j=0$, $j=j_0$, and $j_0=0$.  When $j=0$ \eqref{exp} becomes
\begin{equation}\label{common}
 2(n-1)q j_0 (\frac{1}{p} - \frac{1}{\ptil}) \leq - \eps q j_0,
\end{equation}
which is true for some $\eps > 0$ since $\ptil < p$.  
When $j=j_0$ \eqref{exp} becomes
$$ 2(n-1)q (\frac{1}{p} - 1 + \frac{n+1}{2(n-1)q}) \leq 0,$$
which holds since $p^\prime \leq \frac{n-1}{n+1} 2q$.
Finally, when $j_0 = 0$
\eqref{exp} becomes
\begin{equation}\label{common2}
(2n - 2(n-1)q)j \leq -\eps q j,
\end{equation}
which holds for some $\eps > 0$ since $2q > \frac{2n}{n-1}$.

It remains to treat the case $\ptil > 2q$.  By repeating the
above procedure but with \eqref{xr-big} replaced by \eqref{xr-small},
we see that \eqref{major}
reduces to
\begin{equation}\label{exp2}
\begin{split}
\frac{2(n-1)q}{p}j_0  
-\frac{2(n-1)q}{\ptil^\prime}j& + (n+1)j
-(n-1)\frac{2q}{\ptil}j_0\\
&+(n-1)j (1 - \frac{2q}{\ptil})
\leq -\eps q |j - j_0|.
\end{split}
\end{equation}
Since the left-hand side is completely linear it suffices to
verify this when $j=0$ and when $j_0=0$.  But in these two cases
\eqref{exp2} reduces \eqref{common}, \eqref{common2} as before,
and so the argument proceeds as in the previous case.
\endprf

The fact that this theorem requires knowledge of $R^*(p \times p \to q)$
for all elliptic phase functions
is a defect of the argument.  When
restricted to the quadratic phase $\Phi(x) = \frac{1}{2}|x|^2$ however,
no other phase functions are required in the proof, due to the
algebraic properties of $\Phi$.  The quadratic
phase is the simplest of all the elliptic phases; indeed, a parabolic
scaling and limiting argument shows that any sharp restriction theorem
for an elliptic phase implies the corresponding estimate for the
quadratic phase.  (See \cite{tao:boch-rest}).

\subsection{Necessity of \eqref{c0}-\eqref{c2}}\label{nec}
\

In this section we sketch the proof of the
assertion that the conditions in Conjecture \ref{conj} are
necessary.  For simplicity we take $\Phi$ to be a graphing function
for a small portion of a sphere;
one can easily
modify the arguments below for more general phases.  The estimate
$R^*(p \times p \to q)$ can then be rewritten as
\begin{equation}\label{spherical}
 \| \widehat{fd\sigma} \widehat{gd\sigma} \|_q \lesssim \|f\|_p \|g\|_p,
 \end{equation}
where $d\sigma$ is surface measure on the unit sphere $S^{n-1}$, and
$f$ and $g$ are functions on fixed disjoint caps $C_1$, $C_2$
in $S^{n-1}$ whose size and separation are comparable to a small quantity
$\epsilon = \epsilon_n$.

To prove \eqref{c0}, we take $f(w) = 1$ on $C_1$, and $g(w) = e^{-2\pi i 
x_0 \cdot w}$ on $C_2$, where $x_0 \in \R^n$ is a point to
be determined later.  From standard stationary phase estimates, we see that
for any $R \gg 1$ one can find a cube $Q$ of sidelength $R$ such that
$|\widehat{fd\sigma}(x)| \sim R^{-\frac{n-1}{2}}$ on $Q$.  By choosing
$x_0$ appropriately, one can also arrange matters so that
$|\widehat{gd\sigma}(x)| \sim R^{-\frac{n-1}{2}}$ on the same cube $Q$.
By inserting these estimates into \eqref{spherical} one obtains
$$ R^{-\frac{n-1}{2}}  R^{-\frac{n-1}{2}} |Q|^{1/q} \lesssim 1.$$
If one now uses the fact that $|Q| \sim R^n$ and takes $R \to \infty$
the condition \eqref{c0} follows.

The proof of the necessity of \eqref{c1} and \eqref{c2} is based on
modifications of the standard Knapp example.  We note in passing that
without modification the Knapp example only gives the weaker condition
$$ \frac{n}{2q} + \frac{n-1}{p} \leq n-1.$$

To prove \eqref{c1}, we will take $f$ and $g$ to be ``squashed caps''\footnote{
This example was discovered independently by the authors and Sergiu Klainerman.}.  We factor $\R^n$ as $\R^2 \times \R^{n-2}$, and use $S^1$ to denote
the great circle $S^1 = S^{n-1} \cap (\R^2 \times \{0\})$.  We may assume that
$S^1$ intersects $C_1$ and $C_2$.
Fix $0 < \delta \ll 1$.  We take $f$ and $g$ to
be the characteristic functions of the
sets
$$ C_i \cap (B_2(w_i,\delta^2) \times B_{n-2}(0,\delta)),
\quad i = 1,2$$
respectively, where $B_k(x,R)$ denotes the ball in $\R^k$ of radius $R$
centered at $x$, and $w_i$ are arbitrary elements of $S^1 \cap C_i$ for
$i=1,2$.  Then the Fourier transforms of $fd\sigma$, $gd\sigma$
exhibit essentially no cancellation on the box
\begin{equation}\label{box}
B_2(0,\frac{1}{C\delta^2}) \times B_{n-2}(0, \frac{1}{C\delta}).
\end{equation}
Indeed, we have
$|\widehat{fd\sigma}(x)| \sim |\widehat{gd\sigma}(x)| \sim \delta^n$
on this set.
Inserting this estimate into \eqref{spherical} one obtains
$$ \delta^n \delta^n (\delta^{-n-2})^{1/q} \lesssim 
\delta^{\frac{n}{p}} \delta^{\frac{n}{p}},$$
and by taking $\delta \to 0$ one obtains \eqref{c1}.

The estimate \eqref{c2} is proven by taking $f$ and $g$ to be ``stretched 
caps''.
With the notation as before we take $f$ and $g$ to be the characteristic
functions of 
$$ C_i \cap (\R^2 \times B_{n-2}(0,\delta)),
\quad i = 1,2$$
respectively to begin with, although we will later need to multiply $f$
and $g$ by a phase as in the proof of \eqref{c0}.

When restricted to the slab $\R^2 \times B_{n-2}(0, \frac{1}{C\delta})$,
the functions $\widehat{fd\sigma}, \widehat{gd\sigma}$ behave essentially
like Fourier transforms of measures on $S^1$.  Indeed, a stationary phase
computation shows that
$$ |\widehat{fd\sigma}(x)| \sim \delta^{n-2} |x|^{-\frac{1}{2}}$$
on a large portion of this slab, and similarly for $\widehat{gd\sigma}$.
Thus, multiplying by a phase to
translate $\widehat{fd\sigma}$ and $\widehat{gd\sigma}$ as necessary,
one can arrange matters so that
$$ |\widehat{fd\sigma}(x)| \sim
|\widehat{gd\sigma}(x)| \sim \delta^{n-2} (\delta^{-2})^{-\frac{1}{2}}
= \delta^{n-1}$$
on the box \eqref{box}.   Inserting this into \eqref{spherical} one obtains
$$ \delta^{n-1} \delta^{n-1} (\delta^{-n-2})^{1/q} \lesssim 
\delta^{\frac{n-2}{p}} \delta^{\frac{n-2}{p}},$$
and \eqref{c2} follows by taking $\delta \to 0$.
\endprf

Unlike the situation with the disc multiplier problem \cite{feff:ball},
it appears that the Besicovitch set construction
does not give any further restrictions on $p$, $q$.  Indeed, for $n=2$
the conditions \eqref{c0}-\eqref{c2} are sufficient as well as necessary.

\subsection{Proof of Theorem \ref{modest}}
\

Our argument will be a routine modification of the one in \cite{vargas:2}.

The necessity of the condition on $p$ follows from Section \ref{nec}, so
we will only show the sufficiency of this condition.  By H\"older's inequality it suffices to 
show that
$$ R^*(\frac{4n}{3n-2} \times \frac{4n}{3n-2} \to 2).$$
By symmetry and interpolation this will follow from
$$ R^*(2 \times \frac{n}{n-1} \to 2).$$
It suffices to show that
\begin{equation}\label{interp}
\int \Re^* f_1(x) \Re^*g_1(x) \overline{\Re^* f_2(x)}
\overline{\Re^* g_2(x)}\ dx \lesssim \|f_1\|_{1} \|g_1\|_{\frac{n}{n-1}} 
\|f_2\|_{\infty} \|g_2\|_{\frac{n}{n-1}}
\end{equation}
for all $f_1$, $f_2$, $g_1$, $g_2$, supported on $Q_1$, $Q_1$, $Q_2$, $Q_2$
respectively, where $Q_1$ and $Q_2$ are $O(1)$-separated cubes.
Indeed, by applying the symmetry
$f_1 \leftrightarrow f_2$, $g_1
\leftrightarrow g_2$ to \eqref{interp} and applying multi-linear interpolation
one obtains
$$ \int \Re^* f_1(x) \Re^*g_1(x) \overline{\Re^* f_2(x)}
\overline{\Re^* g_2(x)}\ dx \lesssim \|f_1\|_{2} \|g_1\|_{\frac{n}{n-1}} \|f_2\|_{2} \|g_2\|_{\frac{n}{n-1}},$$
and the desired estimate follows from substituting $f_1=f_2=f$, $g_1=g_2=g$.

It remains to prove \eqref{interp}.  By Plancherel's theorem the left-hand
side can be written as
$$ \int f_1(x) g_1(y) \overline{f_2(z)} \overline{g_2(w)} \delta(\Phi(x)
+\Phi(y) -\Phi(z)-\Phi(w))\delta(x+y-z-w)\ dx dy dz dw,$$
where $\delta$ is the Dirac distribution.  From the positivity of the
kernel in the above expression we may reduce \eqref{interp} to
\begin{align*}
\int f_1(x) g_1(y) f_2(z) g_2(w) &\chi_1(x) \chi_2(y) \chi_1(z) \chi_2(w)\\
&\delta(\Phi(x)
+\Phi(y) -\Phi(z)-\Phi(w)) \delta(x+y-z-w)\ dx dy dz dw \\
&\lesssim
\|f_1\|_{1} \|g_1\|_{\frac{n}{n-1}}
\|f_2\|_{\infty} \|g_2\|_{\frac{n}{n-1}}
\end{align*}
for arbitrary functions $f_1$,$f_2$,$g_1$,$g_2$ on $\R^{n-1}$, where
$\chi_1$ and $\chi_2$ are smooth cutoff 
functions adapted to (a slight thickening of)
$Q_1$ and $Q_2$ respectively.
Since $f_1$ and $f_2$ are controlled in $L^1$ and $L^\infty$ respectively,
we may assume that $f_1(x) = \delta(x-x_0)$ and $f_2 \equiv 1$ 
for some $x_0$; we may take $x_0 = 0$ by translating $\Phi$ and 
subtracting off a harmless affine factor.  In particular, we may
assume that $0$ is in (a slight thickening of) $Q_1$.  The estimate 
\eqref{interp} thus reduces to 
\begin{align*}
\int g_1(y) g_2(w) \chi_2(y) \chi_1(y-w) \chi_2(w)
&\delta(\Phi(y) - \Phi(y-w) - \Phi(w))\ dy dw\\
&\lesssim \| g_1\|_{\frac{n}{n-1}} \|g_2\|_{\frac{n}{n-1}},
\end{align*}
which by duality becomes
\begin{equation}\label{standard}
\| T g\|_n
\lesssim \| g\|_{\frac{n}{n-1}},
\end{equation}
where $T$ is the averaging operator
$$
Tg(y) = \int g(w) \chi_2(y) \chi_1(y-w) \chi_2(w)
\delta(\Phi(y) - \Phi(y-w) - \Phi(w))\ dw.
$$
It is well known (see below) that the estimate \eqref{standard} will
hold if
the defining function $\phi(y,w) = \Phi(y) - \Phi(y-w) - \Phi(w)$
satisfies the rotational curvature condition
\begin{equation}\label{rot}
\left|\det \left(\begin{array}{ll}
\phi& \phi_{y}\\
\phi_{w}& \phi_{y w}\\
\end{array}\right)\right| > 0 \quad \hbox{when } \phi = 0
\end{equation}
uniformly on the support of $\chi_2(y)\chi_1(y-w)\chi_2(w)$.

However, since $\Phi$ is elliptic, we have the estimates
$$ 
\Phi_{y_i y_j} = \delta_{ij}
+ O(\epsilon_0),\quad
\Phi_{y_i}(y) = y_i + O(\epsilon_0 |y|), \quad 
\Phi(y) = \frac{1}{2}|y|^2 + O(\epsilon_0 |y|^2)$$
where $\delta_{ij}$ is the Kronecker delta.  Inserting these estimates
into the definition of $\phi$, one can estimate the above determinant
as
$$ \left|\det \left(\begin{array}{ll}
\phi& \phi_{y}\\
\phi_{w}& \phi_{y w}\\
\end{array}\right)\right| = |w|^2 + O(\epsilon_0 (|y|^2 + |w|^2)) \quad
\hbox{ when } \phi = 0.
$$
However, from the support of $\chi_2(y) \chi_2(w)$ and the assumption
that $0$ is in a thickening of $Q_1$ we see that $|w|, |y| 
\sim 1$.  Thus \eqref{rot} follows, if 
$\epsilon_0$ is sufficiently small.  This finishes the proof.
\endprf

The above proof shows that there exist asymmetrical
bilinear restriction theorems in addition to the symmetrical ones.
In particular, one may conjecture that
$$ R^*(\frac{n+2}{n} \times \frac{n+2}{2} \to \frac{n+2}{n}),$$
which is a strengthening of \eqref{kl}.  Non-symmetrical
versions of the counterexamples in the previous section show that this
conjecture is best possible.

For $n\leq 3$ Theorem \ref{modest} is an improvement on the classical
Tomas-Stein theorem.  However for $n > 3$ the two estimates are not
directly comparable.  Because of this, we have no significant improvements
to Wolff's restriction theorem \cite{wolff:kakeya} in four and higher 
dimensions.

For completeness we sketch a proof of the following standard fact which
was used in the above proof.

\begin{lemma}  If $\phi$ satisfies the rotational curvature condition
\eqref{rot} on the support of a cutoff function $\psi(y,w)$, then the operator
$$ Tf(y) = \int_{\R^{n-1}} f(w) \psi(y,w) \delta(\phi(y,w))\ dw$$
obeys \eqref{standard}.
\end{lemma}

\begin{proof}
We imbed this operator in the analytic family $T_\zeta$ defined by
$$ T_\zeta(y) = \int f(w) \psi(y,w)
a_\zeta(\phi(y,w))\ dw,$$
where $a_\zeta$ is defined for $Re(\zeta) > 0$ by
$$ a_\zeta(t) = e^{\zeta^2} \frac{t_+^{\zeta-1}}{\Gamma(\zeta)} \varphi(t)$$
and $\varphi$ is a cutoff function adapted to $[-\eps,\eps]$ for some
small $\eps > 0$; for $Re(\zeta)\leq 0$
$a_\zeta$ (and thus $T_\zeta$) is defined by analytic continuation.  Since
$T = T_0$, \eqref{standard} will follow from complex interpolation
between the estimates
\begin{align*}
\| T_{1 + it} f \|_\infty &\lesssim \|f\|_1\\
\| T_{-\frac{n-2}{2}+it} f \|_2 &\lesssim \|f\|_2
\end{align*}
for all real $t$ and some fixed $N>0$.
The former estimate follows immediately from the observation that the
kernel of $T_{1+it}$ is uniformly bounded in $t$ (indeed, the $e^{\zeta^2}$
term makes it rapidly decreasing in $t$).  To
prove the latter estimate, it suffices to show that $T_{-\frac{n-2}{2}+it}$
is a Fourier integral operator of order $0$ uniformly 
in $t$ (see e.g. \cite{horm:fio}).
Accordingly, we write $T_{-\frac{n-2}{2}+it}$ as
$$ T_{-\frac{n-2}{2}+it} f(y) =
\int e^{2\pi i \phi(y,z)\cdot \xi} f(w)
\psi(y,w) \hat a_{-\frac{n-2}{2}+it}(\xi) \ dw d\xi$$
where $\xi$ ranges over $\R$.  From the rotational curvature hypothesis
\eqref{rot} we see that the phase is non-degenerate in the sense of
\cite{horm:fio}.  Since the amplitude is a symbol of order $-\frac{n-2}{2}$,
$y$, $w$ range over a $n-1$ dimensional space, and $\xi$ ranges over
a $1$-dimensional space, the reduction-of-variables
theorem (see e.g. \cite{horm:fio})
states that $T_{-\frac{n-2}{2}+it}$ will be a Fourier integral operator
of order $0$, as desired; the uniformity in $t$ follows from the rapid
decrease of $a_{-\frac{n-2}{2}+it}$ with respect to $t$, caused
by the $e^{\zeta^2}$ factor.
\end{proof}

\section{Bilinear Kakeya estimates}\label{b-kakeya}

We now begin the second part of this paper, in which 
we give analogues of the previous results for the Kakeya operator.

Throughout this section
$0 < \delta \ll 1$ will be a small parameter, and
we will use $A \lesssim B$ to denote the estimate $A \leq C_\eps \delta^{-\eps} B$ for all $\eps > 0$, otherwise we write $A \gg B$.  We say that a quantity $A$ has \emph{logarithmic size} if $1 \lesssim |A| \lesssim 1$, while
we say it has \emph{polynomial size}
if $\delta^C \lesssim |A| \lesssim \delta^{-C}$ for some 
constant $C$. 
Finally, all functions and quantities in this section are assumed
to be non-negative.

Let $\E$ be a $\delta$-net of the unit cube $Q$ in $\R^{n-1}$.  We give
two measures on $\E$, counting measure $di$ and normalized
counting measure $d\omega= \delta^{n-1} di$.
For $\omega, i \in \E$, define the $\delta \times 1$ tube $T_{\omega}^i$ by
$$ T_\omega^i = \{(\underline{y},y_n) \in \R^n: |y_n| \leq 1, 
|\underline{y} - y_n \omega - i| \leq \delta\};$$
we will call $\omega$ and $i$ the direction and base of $T_\omega^i$
respectively.
Note that for fixed $\omega$ the tubes $T_\omega^i$
essentially form a partition of the unit ball $B(0,1)$.  This discretization
is not essential to the statements and estimates, but it allows for some
technical simplification to the argument.

For any function $f$ on $\R^n$, define the discretized x-ray 
transform $Xf = X_\delta f$ on $\E \times \E$ by
$$ X f(\omega,i) = \delta^{1-n} \int_{T_\omega^i} f(x)\ dx.$$
For $1\leq p,q \leq \infty$, let $K(p \to q)$ denote the estimate
$$ \| X f\|_{L^q_\omega L^\infty_i} \lesssim \delta^{-\frac{n}{p}+1} 
\|f\|_p,$$
where $\E \times \E$ is understood to be endowed with the
measure $d\omega di$.  By taking $f$ to be the characteristic
function of a $\delta$-ball we see that the factor $\delta^{-\frac{n}{p}+1}$ is
best possible.

The Kakeya conjecture asserts that $K(p \to q)$ holds if and only if
$1 \leq p \leq n$ and $q \leq (n-1)p^\prime$.  In particular, it
is conjectured that $K(n \to n)$ holds.  It is easy to see that
these conditions on $p$, $q$ are necessary.  The conjecture is trivial
for $p=1$; the difficulty is in making $p$ (and to a lesser extent $q$)
as large as possible.  So far the best result
on this conjecture is due to Wolff \cite{wolff:kakeya}, who showed
that 
$$K(\frac{n+2}{2} \to \frac{(n-1)(n+2)}{n}).$$
In particular, for $n=3$ we have $K(\frac{5}{2} \to \frac{10}{3})$.
This estimate is sharp in the sense that the $\frac{10}{3}$ exponent
cannot be raised without decreasing the $\frac{5}{2}$ exponent.

The adjoint estimate 
$$ \| X^* g\|_{\pp} \lesssim \delta^{-\frac{n}{p}+1} \|g\|_{L^\qp_\omega
L^1_i}$$
to $K(p \to q)$ will be denoted $K^*(q^\prime \to p^\prime)$; note that
\begin{equation}\label{xdef}
X^* g(x) = \delta^{1-n}
\int\int g(\omega,i) \chi_{T_\omega^i}(x)\ d\omega di =
\sum_\omega \sum_i
g(\omega,i) \chi_{T_\omega^i}(x).
\end{equation}

Following the philosophy of the previous sections, we define the
bilinear version\footnote{Note that the last exponent
will usually be less than 1.}
$K^*(q^\prime \times q^\prime \to \frac{p^\prime}{2})$
of the above estimate by
$$ \| X^* f X^* g \|_{p^\prime/2}
\lesssim \delta^{-\frac{2n}{p}+2} \|f\|_{L^\qp_\omega L^1_i} 
\|g\|_{L^{q^\prime}_\omega L^1_i}$$
for all 
$f$, $g$ supported on $\E_1 \times \E$, $\E_2 \times \E$, where $\E_1$ and
$\E_2$ are $O(1)$-separated subsets of $\E$.

We have the following analogue of Theorem \ref{imply}, which we 
will prove in Section \ref{kakeya-iff}.  For technical reasons
we will restrict ourselves to the case $p \leq q$, which is the 
case of most interest.  It is likely that one can use factorization
theory and affine invariance to extend these results to the case
$p > q$.

\begin{theorem}\label{kakeya-imply}  Suppose that 
$1 \leq p \leq q \leq (n-1)p^\prime$.  Then $K(p \to q)$ and
$K^*(q^\prime \times q^\prime \to \frac{p^\prime}{2})$ are equivalent.
\end{theorem}

As with the restriction conjecture, it is possible to have
bilinear Kakeya estimates which are outside the range of the usual
Kakeya conjecture.  For instance, one has the easy estimate

\begin{proposition}\label{111}  For any $n \geq 2$ we have
$K^*(1 \times 1 \to 1)$.
\end{proposition}

We defer the simple proof of this proposition to Section \ref{cordoba}.

Interpolating this estimate with the estimate
$$ K^*(\frac{n}{n-1} \times \frac{n}{n-1} \to \frac{n}{2(n-1)}),$$
which by Theorem \ref{kakeya-imply} is the bilinear form of
the Kakeya conjecture, we see that the Kakeya conjecture is equivalent to

\begin{conjecture}\label{kakeya-conj}
If $n \geq 2$ and $1 \leq p,q \leq \infty$, then 
$K^*(q^\prime \times q^\prime \to \frac{p^\prime}{2})$ holds if
and only if
\begin{align}
p &\leq n,\label{k0}\\
\frac{n-2}{q} + \frac{2}{p} &\geq 1.\label{k1}
\end{align}
\end{conjecture}

We will show the necessity of \eqref{k0} and \eqref{k1} in Section
\ref{kakeya-nec}.  These two conditions correspond to \eqref{c0} and
\eqref{c2} respectively; the analogue of \eqref{c1} is the
degenerate condition $q \leq \infty$.

Wolff's theorem \cite{wolff:kakeya} is equivalent to
$$ K^*(\frac{(n-1)(n+2)}{n^2-2} \times \frac{(n-1)(n+2)}{n^2 - 2} \to 
\frac{n+2}{2n});$$
in particular, we have $K^*(\frac{10}{7} \times \frac{10}{7} \to
\frac{5}{6})$ for $n=3$.  In Section \ref{wolff-sec} we 
improve the above estimate to

\begin{theorem}\label{b-wolff}  For all $n \geq 2$ we have
$$ K^*(\frac{n+2}{n+1} \times \frac{n+2}{n+1} \to \frac{n+2}{2n}).$$
\end{theorem}

In particular, we have 
\begin{equation}\label{5-6}
K^*(\frac{5}{4} \times \frac{5}{4} \to
\frac{5}{6})
\end{equation}
in three dimensions.  The result can be thought of as
a bilinear version of the (false) estimate $K(\frac{n+2}{2} \to n+2)$,
and is sharp in the sense that \eqref{k1} is obeyed with equality.

We display the known Kakeya and bilinear Kakeya results in Figure 2.
The trapezoidal region represents the conjectured range of $(p,q)$
for which $K(p \to q)$ should hold, and the pentagonal enlargement
represents the range on which the bilinear version 
$K^*(q^\prime \times q^\prime \to \frac{p^\prime}{2})$ should hold.
By Theorem \ref{kakeya-imply} the two estimates are equivalent in the
triangular region below the dashed line.
The point $1$ is the trivial $L^1 \to L^\infty$ estimate, while
the point $2$ represents the higher-dimensional analogue of Cordoba's
argument (\cite{cordoba:sieve}, \cite{borg:kakeya}), while $3$ is
the ``bush'' argument as given by Bourgain \cite{borg:kakeya}
(see also \cite{drury:xray}, \cite{cordoba:sieve}).  Proposition \ref{111},
the bilinear improvement to Cordoba's argument, is the point $4$.
The point $5$ is Bourgain's Kakeya maximal theorem \cite{borg:kakeya}
(see also \cite{schlag:kakeya}), while $6$ is Wolff's theorem 
\cite{wolff:kakeya}, which we improve in Theorem \ref{b-wolff} to the point
$7$.  The region to the right of the dotted line thus represents the best
results known to date (excepting the results in \cite{wolff:xray},
which are not directly representable on this figure).

\begin{figure}[htbp] \centering
\ \psfig{figure=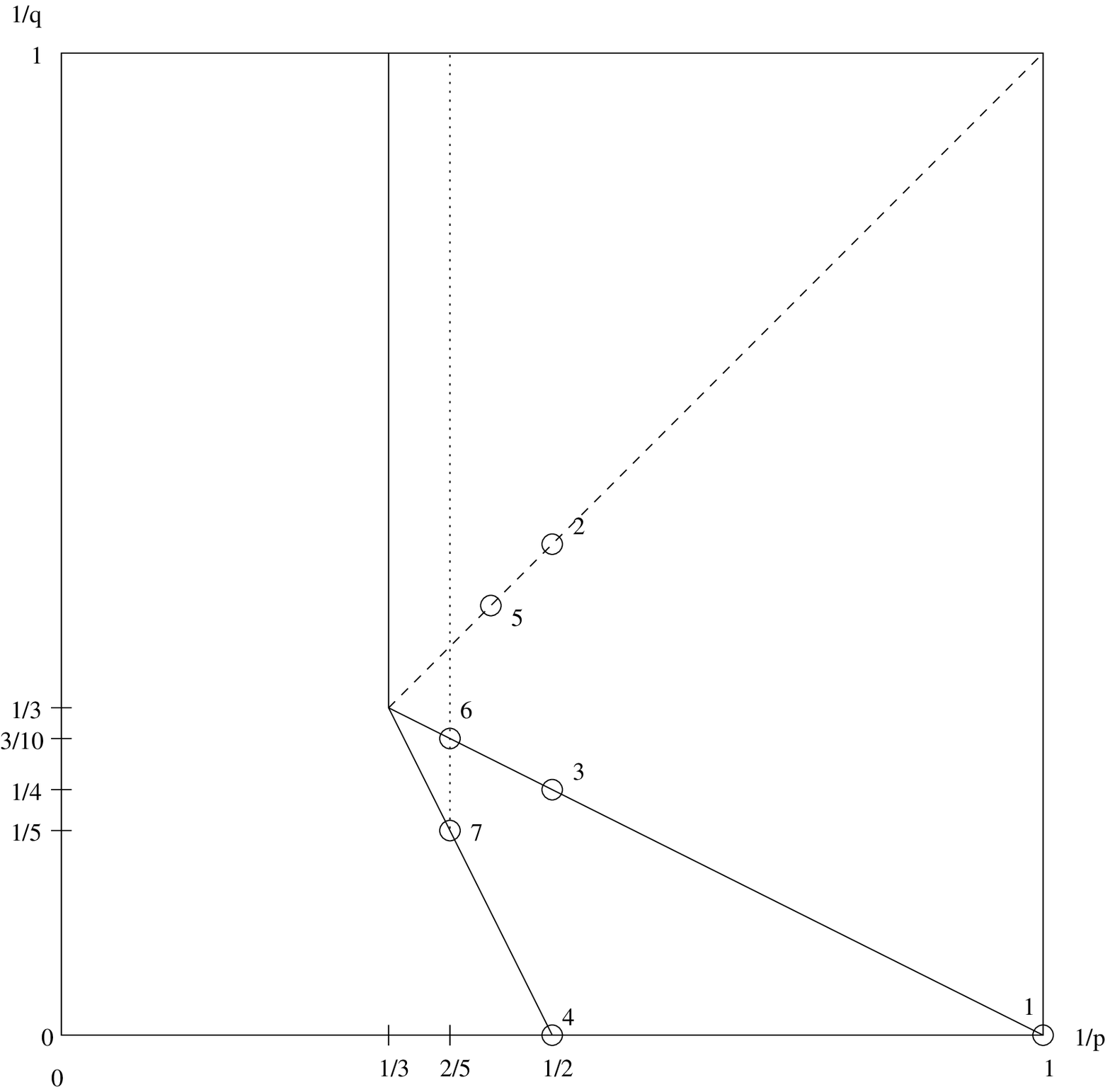,height=3in,width=3.6in}
\caption{Status of $K(p \to q)$ and $K^*(q^\prime \times q^\prime \to \frac{p^\prime}{2}$) for $n=3$.  }
        \label{fig:kakeya}
	\end{figure}

\subsection{Proof of Proposition \ref{111}}\label{cordoba}
\

We will need the following geometric observation of Cordoba:

\begin{lemma}\label{cord-geom}  For any $\omega_1, \omega_2, i_1, i_2 \in \E$,
One has
$$ \langle \chi_{T_{\omega_1}^{i_1}},\chi_{T_{\omega_2}^{i_2}}
\rangle = |T_{\omega_1}^{i_1} \cap T_{\omega_2}^{i_2}| \lesssim
\frac{\delta^{n}}{|\omega_1-\omega_2| + \delta}.$$
\end{lemma}

In particular, if $\omega_1$ and $\omega_2$ have unit separation,
then the intersection between the two tubes has measure at most $\delta^n$.
We leave the easy proof of this lemma to the reader.

From this observation we easily see that
\begin{align*}
\| X^* f X^* g \|_1 &= \langle X^*_\delta f, X^*_\delta g \rangle\\
&= \int\int\int\int \delta^{1-n} \delta^{1-n}
\langle \chi_{T_{\omega_1}^{i_1}},\chi_{T_{\omega_2}^{i_2}}
\rangle f(\omega_1, i_1) g(\omega_2,i_2)\ d\omega_1 di_1 d\omega_2 di_2\\
&\lesssim \int\int\int\int \delta^{1-n} \delta^{1-n} \delta^n
f(\omega_1, i_1) g(\omega_2,i_2)\ d\omega_1 di_1 d\omega_2 di_2\\
&= \delta^{2-n} \|f\|_{L^1_\omega L^1_i} \|g\|_{L^1_\omega L^1_i}
\end{align*}
which is $K^*(1 \times 1 \to 1)$, as desired.
\endprf

\subsection{Proof of Theorem \ref{b-wolff}}\label{wolff-sec}
\

Apart from several technical changes, this theorem will be proven
using the geometric and combinatorial arguments 
of Wolff 
\cite{wolff:kakeya}, namely Cordoba's observation (Lemma \ref{cord-geom})
and the ``brush'' argument.  The bilinear setting allows for some
simplification since the average angular separation $\sigma$ between 
two tubes, as defined in \cite{wolff:kakeya}, may be (heuristically
at least) taken to be $1$.  In fact, this result informally follows
by setting $\sigma=1$ in Lemma 2.1 of \cite{wolff:xray}, and removing
the ``two ends'' condition as in that paper.  We will 
also take advantage of
some simplifications noted by later authors
(notably \cite{schlag:bourgain, schlag:geom-borg,
schlag:kakeya}, \cite{sogge:kakeya}, \cite{wolff:xray}).  Of course, due to the fact that 
we are in a bilinearized
adjoint setting, there are some technical difficulties, most notably
defining
the analogue of the quantity $\lambda$ in \cite{wolff:kakeya}.  
Also, since the target space
$L^{(n+2)/2n}$ is not a Banach space, certain reductions and techniques
(e.g. duality, elimination of the $i_1$, $i_2$ variables) become
unavailable.  In particular, the Lebesgue space approach of \cite{katz:kakeya}
becomes technically very difficult, and we will use
restricted weak-type methods
instead.  In other words, we will use the pigeon-hole principle 
to reduce as many functions as possible to characteristic functions.

We first make the trivial observation that since $X$ is discretized,
the operator boundedness of $X$ on Lebesgue spaces
is automatic with some large power
of $\delta^{-1}$; the issue is to control the dependence on $\delta$
efficiently.

Let us normalize $f$ and $g$ so that
$$
\|f\|_{L^{\frac{n+2}{n+1}}_\omega L^1_i} = 
\|g\|_{L^{\frac{n+2}{n+1}}_\omega L^1_i} = 1.
$$
We have to show that
$$ \| X^* f X^* g\|_{\frac{n+2}{2n}} \lesssim \delta^{-2\frac{n-2}{n+2}}.$$
It will suffice to show the weak-type bound
\begin{equation}\label{weak-1}
 |\{ X^* f X^* g \gtrsim \alpha \}| \lesssim
\alpha^{-\frac{n+2}{2n}} \delta^{-\frac{n-2}{n}}
\end{equation}
for all $\alpha > 0$,
since the strong-type estimate can be recovered (with only a logarithmic
loss) by integrating this over all $\alpha$ of polynomial size; the 
contribution of $\alpha \gg \delta^{-C}$ or $\alpha \ll \delta^C$ can be easily 
controlled using trivial estimates.

We now make the assumption that
there exists sets $\Omega_j \subset \E_j$ of
cardinality $M_j > 0$ for $j=1,2$ such that
\begin{equation}\label{norm}
\| f(\omega,\cdot)\|_{L^1_i} = 
(M_1 \delta^{n-1})^{-\frac{n+1}{n+2}} \chi_{\Omega_1}(\omega),\quad
\| g(\omega,\cdot)\|_{L^1_i} = 
(M_2 \delta^{n-1})^{-\frac{n+1}{n+2}} \chi_{\Omega_2}(\omega).
\end{equation}
This assumption is justified as any $L^{(n+2)/(n+1)}$-normalized
$f$, $g$ can be majorized
by a sum of at most logarithmically many functions of this type.
We may assume that the $M_j$ have polynomial size.

From the pigeon-hole principle \eqref{weak-1} will follow from
the estimate
\begin{equation}\label{goal}
|E| \lesssim (\alpha_1 \alpha_2)^{-\frac{n+2}{2n}} \delta^{-\frac{n-2}{n}}
\end{equation}
where $E$ is any set such that
\begin{equation}\label{lower}
X^* f \geq \alpha_1, X^* g \geq \alpha_2 \hbox{ on } E,
\end{equation}
and $\alpha_1, \alpha_2 > 0$ are arbitrary.  We may assume that
$|E|$, $\alpha_1, \alpha_2$ have polynomial size, since this estimate is easily
obtainable (with a large gain) otherwise.

The $\alpha_j$, $j=1,2$ represent a normalized multiplicity of the tubes in
the supports of $f$ and $g$; roughly speaking, they are related to
the quantity $N$ defined in \cite{wolff:kakeya} by the informal
relationship
$$ \alpha_j \approx (M_j \delta^{n-1})^{-\frac{n+1}{n+2}} N_j.$$

Define the quantity $A$ by
\begin{equation}\label{extreme}
|E| = A (\alpha_1 \alpha_2)^{-\frac{n+2}{2n}} \delta^{-\frac{n-2}{n}};
\end{equation}
We have to show that $A \lesssim 1$.  We may assume without
loss of generality that $A$ is essentially minimal in the sense that
\begin{equation}\label{extreme-2}
|E| \lesssim A (\alpha_1 \alpha_2)^{-\frac{n+2}{2n}} \delta^{-\frac{n-2}{n}}
\end{equation}
for all $\alpha_1$, $\alpha_2$, $E$, $f$, $g$ which obey
\eqref{norm}, \eqref{lower}.

To copy Wolff's argument in \cite{wolff:kakeya} we will
need some control on the quantity $|E \cap 
T_{\omega}^{i}|$ if $T_{\omega}^{i}$ is a ``typical''
tube in a direction in $\E_1$.  (In \cite{wolff:kakeya} such control
is automatic as one is not working in the adjoint setting).  
From \eqref{lower} and \eqref{xdef} we have
the pointwise estimate
$$ 
\delta^{1-n} \int\int f(\omega,i) \chi_{T_\omega^i}(x)\ d\omega di
\geq \alpha_1 \chi_E(x).$$
Integrating this on $E$ we obtain
\begin{equation}\label{av}
 \int\int f(\omega,i) |T_\omega^i \cap E|\ d\omega di \geq
 \lambda_1 \delta^{n-1} (M_1 \delta^{n-1})^{\frac{1}{n+2}},
\end{equation}
where $\lambda_j$ is defined for $j=1,2$ by
\begin{equation}\label{lambda}
\lambda_j = \frac{\alpha_j |E|}{(M_j \delta^{n-1})^{\frac{1}{n+2}}}.
\end{equation}
From our assumptions we see that the $\lambda_j$ are of polynomial
size.

The $\lambda_j$ are the analogues of the quantity $\lambda$ in 
\cite{wolff:kakeya}.  Indeed, from \eqref{norm} and \eqref{av} we expect 
$|T_\omega^i \cap E| \sim \lambda_1 \delta^{n-1}$ on the average.
In fact, because we are considering only an extremal configuration,
a more precise statement is possible.  We say that a tube
$T_\omega^i$ is \emph{good} if $|T_\omega^i \cap E| \geq \frac{1}{4} \lambda_1
\delta^{n-1}$. 
Let $G$ be the set of all $(\omega, i)$ in the support of $f$
associated to good tubes.  The following improvement of \eqref{av}
states that most tubes are good.

\begin{proposition}\label{ext}  We have
$$ \int\int_{G} f(\omega,i)
\ d\omega di\sim (M_1 \delta^{n-1})^{\frac{1}{n+2}}.$$
In particular, we have that $G$ is non-empty, so that $\lambda_1 \lesssim 1$.
\end{proposition}

\begin{proof}  The upper bound follows immediately from \eqref{av},
so it suffices to show the lower bound.  
Let $c > 0$ be a small number of logarithmic size to be
chosen later.  If the lower bound failed, then we would have
$$ \int\int_{G} f(\omega,i)
\ d\omega di \leq c (M_1 \delta^{n-1})^{\frac{1}{n+2}}.$$
The idea is to then replace $f$ by $\tilde f = f\chi_G$, 
and contradict the extremality of $A$ in \eqref{extreme-2}.

Of course, we must modify $\tilde f$ further, as well as $E$, $\alpha_1$ and 
$M_1$, in order to retain \eqref{norm} and \eqref{lower}. 
We replace
 $E$ by
$$ \tilde E = \{ x \in E: X^*(f - \tilde f) < \frac{1}{2} \alpha_1\};$$
note that $\tilde f$ obeys \eqref{lower} if $E$ is replaced by
$\tilde E$ and $\alpha_1$ is replaced by $\frac{1}{2} \alpha_1$.

The next step is to show that $\tilde E$ is comparable to $E$ in size.
From the definition of $\tilde E$ we see that
$$ \int_E X^*(f - \tilde f) \geq \frac{1}{2} \alpha_1 |E \backslash \tilde E|.$$
However, we have from \eqref{xdef} and the definition of $\tilde f$ that
$$ \int_E X^*(f - \tilde f) = \delta^{1-n} \int\int_{G^c}
f(\omega,i) |T_\omega^i \cap E|
\ d\omega di.$$
Combining the two estimates and using the definition of $G$ we obtain
$$ \frac{1}{2} \alpha_1 |E \backslash \tilde E|
\leq \delta^{1-n} \int\int  f(\omega,i)  \frac{1}{4}\lambda_1 \delta^{n-1} \ d\omega di.$$
Using \eqref{norm} and \eqref{lambda} this simplifies to
$$ |E \backslash \tilde E| \leq \frac{1}{2} |E|,$$
so that $|\tilde E| \sim |E|$ as desired.

We now have to modify $\tilde f$, $\alpha_1$,
$\tilde E$, and $M_1$ further
so that \eqref{norm} is restored.
From hypothesis we have
$$ \int \|\tilde f(\omega,\cdot)\|_{L^1_i}\ d\omega
= \int\int \tilde f(\omega,i)
\ d\omega di < c (M_1 \delta^{n-1})^{\frac{1}{n+2}}.$$
However, from \eqref{norm} we have
$$ \| \tilde f(\omega,\cdot)\|_{L^1_i} \leq 
(M_1 \delta^{n-1})^{-\frac{n+1}{n+2}}$$
Thus by H\"older's inequality this implies that
$$ \| \tilde f\|_{L^{\frac{n+2}{n+1}}_\omega L^1_i} \leq
c^{\frac{n+1}{n+2}}.$$
Thus, as before, we can find a logarithmic number of functions $\tilde f_k$
which each obey \eqref{norm} for some $M_1^k$, and such that
$$ \tilde f \lesssim c^{\frac{n+1}{n+2}} \sum_k \tilde f_k.$$
This implies that
$$ \sum_k X^* \tilde f_k \gtrsim c^{-\frac{n+1}{n+2}} \alpha_1$$
on $\tilde E$.  Thus, by reducing $\tilde E$ by a logarithmic factor
one can find a $k$ such that
$$ X^* \tilde f_k \gtrsim c^{-\frac{n+1}{n+2}} \alpha_1$$
on the reduced set (which we will still call $\tilde E$).  

Thus \eqref{lower} is satisfied with $f$ replaced by $\tilde f_k$,
$E$ replaced by $\tilde E$, and
$\alpha_1$ replaced by
$\tilde \alpha_1 \sim c^{-\frac{n+1}{n+2}} \alpha_1$.  But from the 
definition of $A$ this implies that
$$
|\tilde E| \lesssim A (\tilde \alpha_1 
\alpha_2)^{-\frac{n+2}{2n}} \delta^{-\frac{n-2}{n}}.$$
Comparing this with \eqref{extreme-2} and our estimates for $\tilde E$ and
$\tilde \alpha$ we thus obtain a contradiction, if $c$ is sufficiently
small.
\end{proof}

From the above proposition, the definition of $G$ and the identity
$$ \int_E X^*(f\chi_G) = \delta^{1-n} \int\int_{G} f(\omega,i)|T_\omega^i \cap E|\ 
d\omega di$$
we obtain
$$ \int_E X^*(f\chi_G) \gtrsim (M_1 \delta^{n-1})^{\frac{1}{n+2}}
\lambda_1.$$
From \eqref{lambda} this becomes
$$  \int_E X^*(f\chi_G) \gtrsim \alpha_1 |E|.$$
From \eqref{norm} we thus have
$$  \int X^*(f\chi_G) X^* g \gtrsim \alpha_1 \alpha_2 |E|.$$
Expanding out $X^* g$ using \eqref{xdef} this becomes
$$ \delta^{1-n} \int\int g(\omega,i) (\int_{T_\omega^i} X^*(f\chi_G))\ 
d\omega di \gtrsim \alpha_1 \alpha_2 |E|.$$
On the other hand, from \eqref{norm} we have
$$ \int\int g(\omega,i)\ d\omega di = (M_2 \delta^{n-1})^{\frac{1}{n+2}}.$$
Thus there must exist $(\omega_0, i_0)$ in the support of $g$ such that
\begin{equation}\label{good}
\int_{T_{\omega_0}^{i_0}} X^*(f\chi_G) \gtrsim
\frac{\alpha_1 \alpha_2 \delta^{n-1} |E|}{(M_2 \delta^{n-1})^{\frac{1}{n+2}}} =
\alpha_1 \lambda_2 \delta^{n-1}.
\end{equation}
The tube $T_{\omega_0}^{i_0}$ plays the role of the central tube of a
``brush''.  Unlike Wolff's argument in \cite{wolff:kakeya}
(which considered more general
angular separations $\sigma$ than the unit separation), we will be able
to obtain our estimate using only a single brush.  On the other hand,
by utilizing the extremality hypothesis as in Proposition \ref{ext}, one
could certainly obtain a large number of brushes if desired.

By affine invariance
we may take $\omega_0 = i_0 = 0$, so that the central tube is the
vertical tube through the origin.  In particular, $0$ is in $\E_1$,
so every $\omega$ in $\E_2$ has roughly unit separation from the origin.

Let $G_0 \subset G$ be 
the collection of all good $(\omega, i)$ in the support of $f$ such that
$T_\omega^i$ intersects the central tube $T_{0}^{0}$.  Then
expanding out $X^*(f\chi_G)$ in \eqref{good}, we thus obtain
$$ \int\int_{G_0} f(\omega,i) \delta^{1-n} |T_{0}^{0}
\cap T_\omega^i|\ d\omega di \gtrsim \alpha_1 \lambda_2 \delta^{n-1}.$$
From Lemma \ref{cord-geom} we have $|T_{0}^{0}
\cap T_\omega^i| \lesssim \delta^n$, so that
$$ \int\int_{G_0} f(\omega,i)\ d\omega di \gtrsim \alpha_1 \lambda_2 \delta^{n-2}.$$
Let $\Omega_0 \subset \Omega_1$ be the collection of all $\omega$ such that
$(\omega, i)$ is in $G_0$ for at least one $i$.  From \eqref{norm} we
see that
$$ \int\int_{G_0} f(\omega,i)\ d\omega di \lesssim \# \Omega_0 \delta^{n-1} 
(M_1 \delta^{n-1})^{-\frac{n+1}{n+2}},$$
so that
\begin{equation}\label{omega0}
 \# \Omega_0 \gtrsim \frac{\alpha_1 \lambda_2 \delta^{-1}}
{(M_1 \delta^{n-1})^{-\frac{n+1}{n+2}}}.
\end{equation}
For each $\omega \in \Omega_0$ we choose a tube $T_\omega$ from $G_0$
which is in the direction of $\omega$.  These tubes form the ``bristles''
of the brush.  From construction, $|\omega| \sim 1$,
$T_\omega$ intersects $T_{0}^{0}$, and 
$$
 |T_\omega \cap E| \gtrsim \lambda_1 \delta^{n-1}.
$$
As in  Wolff \cite{wolff:kakeya}, we will use \eqref{good-lambda}
to obtain a lower bound on the size of $E$.  More precisely, we will show that
\begin{equation}\label{wolff}
|E| \gtrsim \# \Omega_0 \lambda_1^n \delta^{n-1}.
\end{equation}
Combining this with \eqref{omega0} yields 
$$ |E| \gtrsim \frac{\alpha_1\lambda_1^n \lambda_2 \delta^{n-2}}
{(M_1 \delta^{n-1})^{-\frac{n+1}{n+2}}}.$$
By a completely symmetrical argument one also has
$$ |E| \gtrsim \frac{\alpha_2\lambda_2^n \lambda_1 \delta^{n-2}}
{(M_2 \delta^{n-1})^{-\frac{n+1}{n+2}}}.$$
Multiplying these estimates together one obtains
$$ |E|^2 \gtrsim \frac{\alpha_1 \alpha_2 (\lambda_1 \lambda_2)^{n+1}
\delta^{2n-4}}{(M_1 \delta^{n-1} M_2 \delta^{n-1})^{-\frac{n+1}{n+2}}}.$$
Applying \eqref{lambda} this reduces to
$$ |E|^2 \gtrsim (\alpha_1 \alpha_2)^{n+2} |E|^{2n+2} \delta^{2n-4},$$
which simplifies to \eqref{goal}, as desired.

It remains to prove \eqref{wolff}.  We use
the argument in \cite{wolff:kakeya}; we adopt
the observation in \cite{sogge:kakeya} (see also \cite{katz:kakeya})
that one does not need to utilize the
``two ends'' reduction in \cite{wolff:kakeya} to achieve \eqref{wolff}.

We need some notation.  For all dyadic numbers $\lambda_1 \lesssim \beta
\lesssim 1$ let $\Gamma_\beta$ be the cylindrical region
$$ \Gamma_\beta = \{(\underline{y}, y_n): |\underline{y}| \sim \beta\}.$$
From the properties of $T_\omega$ we see that
$$ \sum_{\lambda_1 \lesssim \beta
\lesssim 1} |T_\omega \cap E \cap \Gamma_\beta| \gtrsim 
 \lambda_1 \delta^{n-1}$$
 for all $\omega \in \Omega_0$.  By the pigeonhole principle,
one can refine $\Omega_0$ by a logarithmic factor so that
\begin{equation}\label{good-lambda}
  |T_\omega \cap E \cap \Gamma_\beta| \gtrsim 
 \lambda_1 \delta^{n-1}
 \end{equation}
 for all $\omega$ in the refined $\Omega_0$, and some $\lambda_1 \lesssim
 \beta \lesssim 1$
 independent of the choice of $\omega$; henceforth this $\beta$
 is considered fixed.

The directions in $\Omega_0$ are $\delta$-separated.  It will be more
convenient to work with a sparser set of directions, so we take
$\tilde \Omega_0$ to be any $\delta/\beta$-net of $\Omega_0$.  From the estimates
$\# \tilde \Omega_0 \gtrsim \beta^{n-1} \# \Omega_0$
and $\beta \gtrsim \lambda_1$
we see that \eqref{wolff} will follow from
\begin{equation}\label{wolff-alt}
|E| \gtrsim \# \tilde \Omega_{0} \frac{\lambda_1^2}{\beta} \delta^{n-1}.
\end{equation}
Let $\Theta$ be a $\delta/\beta$-net of the unit sphere $S^{n-2}$ in
$\R^{n-1}$.  For each $\omega \in \tilde \Omega_0$, we associate
an (essentially unique) element $\theta = \theta_\omega$ of $\Theta$ by
requiring that
$$|\theta - \frac{\omega}{|\omega|}| \lesssim \delta/\beta;$$
recall that $|\omega| \sim 1$ for all $\omega \in \tilde \Omega_0$. 
Furthermore, from elementary geometry and the fact that
$T_\omega$ intersects $T_0^0$ we see that
$T_\omega\cap \Gamma_\beta$ is contained in
the slab $\Pi_\theta$ given by
$$ \Pi_\theta = \{ (\underline{y}, y_n): |\underline{y}| \sim \beta,
|\frac{\underline{y}}{|\underline{y}|} - \theta | \lesssim \delta/\beta.
\}.$$
As the $\Pi_\theta$ are essentially disjoint, \eqref{wolff-alt} will follow
from the estimate
\begin{equation}\label{wolff-2}
|E \cap \Pi_\theta| \gtrsim \# \tilde \Omega_{0,\theta} \frac{\lambda_1^2}{\beta} \delta^{n-1},
\end{equation}
for all $\theta \in \Theta$, where
$$\tilde \Omega_{0,\theta} = \{\omega \in \tilde \Omega_0: \theta_\omega = \theta\}.$$
For the remainder of the argument $\omega$ (and later $\tilde \omega$)
are always assumed to range over $\tilde \Omega_{0,\theta}$.

We now estimate the quantity
$$ Q = \int_{E \cap \Pi_\theta} \sum_{\omega} \chi_{T_\omega}$$
in two different ways.  Firstly, from
the above geometrical considerations and \eqref{good-lambda} 
we have
$$
 |T_\omega \cap E \cap \Pi_\theta| \gtrsim \lambda_1 \delta^{n-1}$$
for all $\omega$ in $\tilde \Omega_{0,\theta}$.
Summing the above
estimate we obtain
\begin{equation}\label{q}
Q
\gtrsim
\# \tilde \Omega_{0,\theta} \lambda_1 \delta^{n-1}.
\end{equation}
We now obtain a different estimate for $Q$.  From the Cauchy-Schwarz
inequality we have
$$ Q \lesssim
|E \cap \Pi_\theta|^{1/2} (\int_{|E \cap \Pi_\theta|} (\sum_{\omega} \chi_{T_\omega})^2)^{1/2}.$$
Squaring both sides and expanding out the integrand into the diagonal
and off-diagonal term, this reduces to
\begin{equation}\label{cord}
\frac{Q^2}{|E \cap \Pi_\theta|} \lesssim
(\int_{|E \cap \Pi_\theta|} \sum_{\omega} \chi_{T_\omega})
+ \sum\sum_{\omega \neq \tilde\omega}
|T_\omega \cap T_{\tilde \omega} \cap E \cap \Pi_\theta|.
\end{equation}
The first term on the right-hand side is just $Q$.  The second term
we may estimate by Lemma \eqref{cord-geom}.  Thus \eqref{cord} becomes
$$ \frac{Q^2}{|E \cap \Pi_\theta|} \lesssim 
Q + \sum\sum_{\omega \neq \tilde\omega} 
\frac{\delta^n}{|\omega - \tilde \omega|}.$$
However, $\omega$, $\tilde \omega$
range over a $\delta/\beta$-separated
set whose elements are within $\delta/\beta$ of the ray $\R^+ \theta$.
Thus for each $\omega$, the number of $\tilde \omega$ such that
$|\omega - \tilde \omega| \sim 2^{-j}$ is at most $1/(2^j\delta/\beta)$,
for any $j$.  Thus the above estimate reduces to
$$  \frac{Q^2}{|E \cap \Pi_\theta|} \lesssim
Q + \sum_{\delta/\beta \lesssim 2^j \lesssim 1} 
\# \tilde\Omega_{0,\theta} \frac{1}{2^j\delta/\beta} 
\frac{\delta^n}{2^{-j}}.$$
Since the number of such $j$ is only logarithmic, we may simplify
the above to
$$ \frac{1}{|E \cap \Pi_\theta|} \lesssim \frac{1}{Q} + \frac{\# \tilde \Omega_{0,\theta} \beta \delta^{n-1}}{Q^2}.$$
Combining this with \eqref{q} and using the hypothesis $\beta \gtrsim \lambda_1$ we obtain
$$ \frac{1}{|E \cap \Pi_\theta|} \lesssim \frac{\beta}{\# \tilde
\Omega_{0,\theta} \lambda_1^2 \delta^{n-1}},$$
which is \eqref{wolff-2}.
This finishes the proof.
\endprf

\subsection{Proof of Theorem \ref{kakeya-imply}}\label{kakeya-iff}
\

The proof will be a reprise of the argument in Theorem \ref{imply}.
The main difference is that the quasi-orthgonality estimate is replaced
by a quasi-triangle inequality, namely Lemma \ref{young} in the Appendix.  
Also the argument is
technically simpler as we allow a logarithmic loss
in the estimates.  The case $p=1$ is trivial, so we will assume $p>1$.

The implication of $K^*(q^\prime \times q^\prime \to \frac{p^\prime}{2})$
from $K(p \to q)$ is immediate from duality and H\"older's inequality.
Now suppose that $K^*(q^\prime \times q^\prime \to \frac{p^\prime}{2})$
holds for some $p,q$ obeying $1 \leq p \leq n$, $1 \leq q \leq (n-1)p^\prime$.
We have to show that $K^*(\qp \to \pp)$ holds.  Since the Kakeya conjecture
is known to hold for $p \leq 2$ (see e.g. \cite{borg:kakeya}, 
\cite{wolff:kakeya}) we may assume that $p>2$.

Let $f$ be an arbitrary function on $\E \times \E$.
We have to show that
\begin{equation}\label{kt1}
 \| X^* f X^* f \|_{\frac{p^\prime}{2}} \lesssim
\delta^{-\frac{2n}{p}+2}
\|f\|_{L^\qp_\omega L^1_i}^2.
\end{equation}

For each integer $j > 0$ such that $\delta \lesssim 2^{-j}$,
we divide $\E$ into $\sim 2^{(n-1)j}$ dyadic ``subcubes'' $\E \cap 
\tau^j_k$ of sidelength $2^{-j}$, and define the notion of closeness
$\tau^j_k \sim \tau^j_{\kp}$
as in Section \ref{iff}.  We partition $X^* f X^* f$ as
$$ X^* f X^* f =
\sum_j \sum_{k,\kp: \tau^j_k \sim \tau^j_{\kp}} \sum_{m,\mp}
X^*(f \chi_{\tau^j_k} \otimes \chi_{\tau^j_m})
X^*(f \chi_{\tau^j_\kp} \otimes \chi_{\tau^j_\mp}).$$
We now observe the geometric fact that the summand in the above expression
is only non-zero when $\tau^j_m$ and $\tau^j_\mp$ are within $O(2^{-j})$
of each other; we will implicitly assume this in our summation.
By inserting the above identity into \eqref{kt1} and applying Lemma \ref{young}
from the Appendix we reduce ourselves to
\begin{equation}\label{k-target}
 \sum_j \sum_{k,\kp : \tau^j_k \sim \tau^j_{\kp}} \sum_{m,\mp}
\|
X^*(f \chi_{\tau^j_k} \otimes \chi_{\tau^j_m})
X^*(f \chi_{\tau^j_\kp} \otimes \chi_{\tau^j_\mp})
\|_{\frac{p^\prime}{2}}^{\frac{p^\prime}{2}}
\lesssim \delta^{\frac{p^\prime}{2}(-\frac{2n}{p}+2)}
\|f\|_{L^\qp_\omega L^1_i}^\pp.
\end{equation}
This will follow from the following analogue of Proposition \ref{rescale}.

\begin{proposition}\label{k-rescale}  We have
\begin{align*} \|
X^*(f \chi_{\tau^j_k} \otimes \chi_{\tau^j_m})
X^*(f \chi_{\tau^j_\kp} \otimes \chi_{\tau^j_\mp})
\|_{\frac{p^\prime}{2}}^{\frac{p^\prime}{2}}
\lesssim &2^{(1 - \frac{(n-1)\pp}{q})j} \delta^{\frac{p^\prime}{2}(-\frac{2n}{p}+2)}\\
&\|f \chi_{\tau^j_k} \otimes \chi_{\tau^j_m}\|_{L^\qp_\omega L^1_i}^{\frac{p^\prime}{2}}
\|f \chi_{\tau^j_\kp} \otimes \chi_{\tau^j_\mp}\|_{L^\qp_\omega L^1_i}^{\frac{p^\prime}{2}}.
\end{align*}
\end{proposition}

\begin{proof}
By an affine
transformation we may take $\tau^j_k$, $\tau^j_m$
to be centered at the origin.

Applying the hypothesis $K^*(q^\prime \times q^\prime \to \frac{p^\prime}{2})$
to tubes of eccentricity $2^j \delta$ we see that
$$ \| X^*_{2^j \delta} f X^*_{2^j \delta} g \|_{p^\prime/2}
 \lesssim (2^j \delta)^{-\frac{2n}{p}+2} \|f\|_{L^\qp_\omega L^1_i} 
 \|g\|_{L^{q^\prime}_\omega L^1_i}$$
For all $f$, $g$ whose $\omega$-supports are on disjoint cubes.  Applying
the rescaling $(\underline{x}, x_n) \to (2^{-j} \underline{x},x_n)$,
$(\omega, i) \to (2^{-j} \omega, 2^{-j} i)$ to this estimate we obtain
$$
\| X^*_\delta f X^*_\delta g \|_{p^\prime/2}
\lesssim 2^{j(n-1)} 2^{j(n-1)} 2^{-\frac{n-1}{p^\prime/2}j}
 (2^j \delta)^{-\frac{2n}{p}+2} 2^{-\frac{n-1}{q}j} \|f\|_{L^\qp_\omega L^1_i}
  2^{-\frac{n-1}{q}j} \|g\|_{L^\qp_\omega L^1_i}$$
whenever $f$ and $g$ are supported on $\tau^j_k \times \tau^j_m$ and 
$\tau^j_\ktil \times \tau^j_\mp$ respectively, and the proposition follows
from substitution and some algebra.
\end{proof}

From this proposition \eqref{k-target} reduces to
$$ \sum_j \sum_{k,\kp: \tau^j_k \sim \tau^j_{\kp}} \sum_{m,\mp}
2^{(1- \frac{(n-1)\pp}{q})j}
\|f \chi_{\tau^j_k} \otimes \chi_{\tau^j_m}\|_{L^\qp_\omega L^1_i}^{\frac{p^\prime}{2}}
\|f \chi_{\tau^j_\kp} \otimes \chi_{\tau^j_\mp}\|_{L^\qp_\omega L^1_i}^{\frac{p^\prime}{2}}
\lesssim
\|f\|_{L^\qp_\omega L^1_i}^\pp.$$
Since there are only logarithmically many $j$'s it suffices to show this
for a fixed $j$.  By polarization it suffices to show that
$$
 \sum_k \sum_m 
2^{(1- \frac{(n-1)\pp}{q})j} \|f \chi_{\tau^j_k} \otimes \chi_{\tau^j_m}\|_{L^\qp_\omega L^1_i}^\pp \lesssim \|f\|_{L^\qp_\omega L^1_i}^\pp,
$$
which we rewrite as
\begin{equation}\label{headache}
 2^{(\frac{1}{\pp} - \frac{n-1}{q})j} \left(\sum_k \sum_m 
\|f_{k,m}\|_{L^\qp_\omega L^1_i}^\pp\right)^{1/\pp} 
\lesssim \|f\|_{L^\qp_\omega L^1_i}
\end{equation}
where $f_{k,m} = f \chi_{\tau^j_k} \otimes \chi_{\tau^j_m}$.  
It suffices to verify this for the case $p=1$ and for the
endpoint $(p,q) = (n,n)$, since the general case $1 \leq p \leq q 
\leq (n-1) p^\prime$ follows by interpolation.  In these two cases
\eqref{headache} becomes
\begin{align}
2^{-\frac{n-1}{q} j} \sup_{k,m} \| f_{k,m}\|_{L^q_\omega L^1_\omega}
&\lesssim \| f\|_{L^q_\omega L^1_\omega}\label{1q}\\
\left(\sum_k \sum_m \|f_{k,m}\|_{L^\np_\omega L^1_i}^\np\right)^{1/\np}
&\lesssim \|f\|_{L^\np_\omega L^1_i}\label{nn}
\end{align}
respectively.  The estimate \eqref{1q} is trivial, while
\eqref{nn} follows from a further interpolation between
the trivial estimates
\begin{align*}
\sup_{k,m}\|f_{k,m}\|_{L^\infty_\omega L^1_i}
&\lesssim \|f\|_{L^\infty_\omega L^1_i}\\
\sum_k \sum_m \|f_{k,m}\|_{L^1_\omega L^1_i}
&\lesssim \|f\|_{L^1_\omega L^1_i}.
\end{align*}
\endprf

We note that if one inserts the result of Theorem \ref{b-wolff} into
the above line of reasoning, then one not only recovers Wolff's 
Kakeya estimate, but also the entropy estimate improvement proven in
Lemma 2.1 of \cite{wolff:xray}.

\subsection{Necessity of \eqref{k0}-\eqref{k1}}\label{kakeya-nec}
\

We now show that the assumptions \eqref{k0} and \eqref{k1} in
Conjecture \ref{kakeya-conj} are necessary. 

To show the necessity of \eqref{k0}, we take 
$$ f(\omega,i) = \chi_{\E_1}(\omega) \delta_{i,0},\quad
 g(\omega,i) = \chi_{\E_2}(\omega) \delta_{i,i_0}$$
where $\delta_{i,j}$ denotes the Kronecker delta function, and $i_0$ is a suitable
point.  A routine calculation using \eqref{xdef}
shows that (if $i_0$ is chosen properly) $X^* f, X^* g \sim 1$ on a ball of radius $\sim 1$.
Inserting this into $K^*(q^\prime \times q^\prime \to \frac{p^\prime}{2})$ yields
$$ 1
\lesssim \delta^{-\frac{2n}{p}+2},$$
and by taking $\delta \to 0$ one obtains \eqref{k0}.

To show the necessity of \eqref{k1}, we adapt the ``stretched
caps'' example used to show \eqref{c2}.  We consider a tube
$T = \R \times B_{n-2}(0,\delta)$ in $\R^{n-1}$,
and take
$$ f(\omega,i) = \chi_{\E_1 \cap T}(\omega) \delta_{i,i_1},\quad
 g(\omega,i) = \chi_{\E_2 \cap T}(\omega) \delta_{i,i_2}.$$
If $i_1$ and $i_2$ are chosen appropriately, then
$X^* f, X^* g$ are both comparable to $1$ on a slab
which looks roughly like $B_2(0,1) \times B_{n-2}(0,\delta)$.
Inserting this into 
$K^*(q^\prime \times q^\prime \to \frac{p^\prime}{2})$ yields
$$ \delta^{\frac{n-2}{p^\prime/2}}
\lesssim 
\delta^{-\frac{2n}{p}+2}
\delta^{\frac{n-2}{\qp}} \delta^{\frac{n-2}{\qp}},$$
and by taking $\delta \to 0$ one obtains \eqref{k1}.

\section{Applications}\label{apply}

In this section we use the bilinear estimates above to prove the following
restriction theorems.

\begin{theorem}\label{main} If $n=3$, then $\R^*(p \to q)$ holds
whenever $p > \frac{170}{77}$ and $q > \frac{34}{9}$.  Furthermore,
$\R^*_s(q)$ holds for all $q > 4 - \frac{5}{27}$.
\end{theorem}

The proof of this theorem will be based on bilinear versions of
certain arguments of Bourgain (\cite{borg:kakeya}, \cite{borg:stein};
see also \cite{vargas:restrict}).  The first step will be to obtain
localized linear and bilinear restriction theorems.

\begin{definition}  If $1 \leq p, q \leq \infty$ and $\alpha \geq 0$, then
we use $\R^*(p \to q, \alpha)$ to denote the estimate
$$ \| \Re^* f \|_{L^q(B_R)} \lesssim R^{\alpha} \|f\|_p,$$
and $\R^*(p \times p \to q,\alpha)$ to denote the estimate
$$ \| \Re^* f \Re^* g \|_{L^q(B_R)} \lesssim R^{\alpha} \|f\|_p \|g\|_p$$
where $f, g, \Re^*$ are as in Section \ref{rest-sec}
and $B_R$ is a ball of radius $R$ in $\R^n$ (the center of $B_R$ 
is irrelevant by translation symmetry).
\end{definition}

It will be convenient to recast these estimates as a restricted
bilinear estimate on the Fourier transform.

\begin{proposition}\label{rest-equiv}
$R^*(p \times p \to q, \alpha)$ is true if and
only if one has
\begin{equation}\label{reform}
 \| \hat f \hat g\|_{L^p(B(x,R))} \lesssim R^\alpha R^{-1/\pp} \|f\|_p R^{-1/\pp}
 \|g\|_p
 \end{equation}
for all $R \gg 1$, $x \in \R^n$ and all functions $f$, $g$ supported on
$A_1^R$, $A_2^R$, where
$$ A_i^R = \{ (\x, \Phi(\x)+t): \x \in Q_i, |t| \lesssim R^{-1} \}.$$
\end{proposition}

\begin{proof}
If $R^*(p \times p \to q, \alpha)$ holds, then \eqref{reform} follows by
translating $\Phi$ by $O(1/R)$ and averaging using 
H\"older's inequality.  Now suppose that \eqref{reform} holds.  To show
$R^*(p \times p \to q, \alpha)$ it suffices to show that
$$ \| (\hat \phi_R \Re^* f) (\hat \phi_R \Re^* g) \|_q \lesssim R^\alpha
\|f\|_p \|g\|_p,$$
where $\phi_R$ is a real radial $L^1$-normalized bump function 
adapted to $B(0,C/R)$, such that $\hat \phi_R$ is non-negative on $B(x,R)$.
But this follows
from \eqref{reform}, Young's inequality, and the identity
$$ \hat \phi_R \Re^* f  = \widehat{\tilde f d\sigma * \phi_R}$$
where $\tilde f(\underline{x},\Phi(\underline{x})) = f(\underline{x})$ is 
the lift of $f$ to the surface
$\{(\underline{x}, \Phi(\underline{x})): \underline{x} \in Q\}$.
\end{proof}

From this proposition and the trivial estimate
$$ \|\hat{f} \hat{g}\|_1 \leq \|\hat{f}\|_2 \|\hat{g}\|_2 = \|f\|_2 \|g\|_2$$
we obtain the bilinear trace lemma
\begin{equation}\label{trace}
\R^*(2 \times 2 \to 1, 1).
\end{equation}
Thus by interpolating this with other estimates (such as \eqref{12-7})
we may obtain estimates of the form $\R^*(p \times p \to q, \alpha)$
with a large value of $\alpha$.  To lower the value of $\alpha$
we will use a bilinear form of an 
argument of Bourgain \cite{borg:kakeya, borg:stein} (see also
\cite{vargas:restrict}):

\begin{lemma}\label{b1} If $2 < p, q < \infty$ and $\alpha > 0$ are such that
$K^*(\frac{p}{2} \times \frac{p}{2} \to \frac{q}{2})$ and $R^*(2 \times 2
\to q, \alpha)$ hold, then $R^*(p \times p \to q, \frac{\alpha}{2} + \eps)$ 
holds for all $\eps > 0$.
\end{lemma}

\begin{proof}
From Proposition \ref{rest-equiv} it suffices to show that
\begin{equation}\label{step-1}
 \| \hat{f} \hat{g} \|_{L^q(B(0,R^2))} \lesssim R^{\alpha+\eps} R^{-2/\pp} 
\|f\|_p R^{-2/\pp} \|g\|_p
\end{equation}
for all $f$, $g$ supported on $A_1^{R^2}$, $A_2^{R^2}$ respectively, 
and all $\eps > 0$.  (The implicit constants will depend on $\eps$).

Let $\phi_R$ be as in Proposition \ref{rest-equiv}, and define
$\phi_R^x(\xi) = e^{-2\pi i x \cdot \xi} \phi_R(\xi)$ for all $x \in \R^n$.  
Then from the
hypothesis $R^*(2 \times 2 \to q, \alpha)$ and Proposition \ref{rest-equiv}
we have
$$ \| \widehat{\phi_R^x} \hat{f} \widehat{\phi_R^x} \hat{g} \|_{L^q(B(x,R))}
\lesssim R^\alpha R^{-1/2} \|f * \phi_R^x\|_2 R^{-1/2} \|g * \phi_R^x\|_2$$
for all $x$.  Averaging this over all $x \in B(0,R^2)$ we obtain
$$ \| \hat{f} \hat{g} \|_{L^q(B(0,R^2))} \lesssim R^\alpha 
( R^{-n} \int_{B(0,R^2)} (R^{-1/2} \|f * \phi_R^x\|_2 R^{-1/2} 
\|g * \phi_R^x\|_2)^q\ dx)^{1/q}.$$
Thus to show \eqref{step-1} it suffices to show that
\begin{equation}\label{step-2}
  R^{-n} \int_{B(0,R^2)} (R^{-1/2} \|f * \phi_R^x\|_2 R^{-1/2} 
\|g * \phi_R^x\|_2)^q\ dx
\lesssim R^\eps (R^{-2/\pp} \|f\|_p R^{-2/\pp} \|g\|_p)^q.
\end{equation}
This will be accomplished by repeated use of the uncertainty principle
and Plancherel's theorem, together with the Kakeya hypothesis.

Let $\E$, $\E_1$, $\E_2$ be as in Section \ref{b-kakeya} with 
$\delta = \frac{1}{R}$.  We partition the annuli $A_i^{R^2}$ into caps
$C_\omega$ for $\omega \in \E_i$, defined by
$$ C_\omega = \{ (\x, x_n) \in A_i^{R^2}: -\nabla \Phi(x) = \omega + O(\frac{1}{R}) \}.$$
From the ellipticity of $\Phi$ and some elementary geometry we see that the 
$C_\omega$ are essentially disks of diameter $1/R$ and thickness $1/R^2$
oriented in the direction $(\omega,1)$,
which form a finitely overlapping cover of $A_i^{R^2}$.  We decompose
$$ f = \sum_{\omega \in \E_1} f_\omega, \quad g = \sum_{\omega \in \E_2} g_\omega$$
where $f_\omega$, $g_\omega$ are adapted restrictions of $f$, $g$ respectively 
to (a suitable dilate of) $C_\omega$.

From the support conditions on $f_\omega$, $g_\omega$ and $\phi_R^x$ we
see that \eqref{step-2} reduces to
\begin{equation}\label{step-3}
\begin{split}
  R^{-n} \int_{B(0,R^2)} (R^{-1/2} 
  (\sum_{\omega \in \E_1} \|f_\omega * \phi_R^x\|_2^2)^{1/2} &
  R^{-1/2}  
  (\sum_{\omega \in \E_2} \|g_\omega * \phi_R^x\|_2^2)^{1/2} 
  )^q\ dx\\
  & \lesssim
  R^\eps (R^{-2/\pp} \|f\|_p R^{-2/\pp} \|g\|_p)^q
\end{split}
  \end{equation}
The function $\widehat{\phi_R^x}$ is rapidly decreasing outside of the
ball $B(x,R)$.  Thus by Plancherel's theorem the left-hand side of
\eqref{step-3} is majorized by
\begin{equation}\label{step-4}
 R^{-n} \int_{B(0,R^2)} (R^{-1/2} 
(\sum_{\omega \in \E_1} \|\widehat{f_\omega}\|_{L^2(B(x,R))}^2)^{1/2} 
 R^{-1/2}
(\sum_{\omega \in \E_2} \|\widehat{g_\omega}\|_{L^2(B(x,R))}^2)^{1/2}
	  )^q\ dx,
\end{equation}
since the portions of $\widehat{\phi_R^x}$ on translates of $B(x,R)$ can
be handled by translation symmetry.

Let $\psi_\omega$ be a Schwarz function which is comparable to $1$ on
$C_\omega$ and rapidly decreasing away from this cap, and whose Fourier
transform satisfies the pointwise estimate
$$ |\hat \psi_\omega(x)| \lesssim  R^{-n-1} \chi_{R^2 \tilde T^\omega_0}(x),$$
where $\tilde T^\omega_0$ is a thickening of $T^\omega_0$, and $R^2
\tilde T^\omega_0 = \{R^2 x: x \in \tilde T^\omega_0\}$.

If we define $\tilde f_\omega = f_\omega / \psi_\omega$, we have
the estimate
$$|\widehat{f_\omega}(x)| = 
|\widehat{\tilde f_\omega} * \widehat{\psi_\omega}(x)|
\lesssim R^{-n-1} \int_{x+R^2 \tilde T_{\omega_0}} |\widehat{\tilde f_\omega}(y)|\ dy.$$
From H\"older's inequality and \eqref{xdef} we thus obtain
\begin{align*}
|\widehat{f_\omega}(x)|^2 &\lesssim
R^{-n-1} \int_{x+R^2 \tilde T_{\omega_0}} |\widehat{\tilde f_\omega}(y)|^2\ dy\\
&\lesssim X^* F_\omega(\frac{x}{R^2})
\end{align*}
where 
$$F_\omega(\omega^\prime,i) = \delta_{\omega,\omega^\prime} 
R^{n-1} \int_{\tilde T^\omega_i} |\widehat{\tilde f_\omega}(R^2 x)|^2\ 
dx$$
and $\delta_{\omega,\omega^\prime}$ is the Kronecker delta.

Since $X^* F_\omega$ is essentially constant on balls of radius $1/R$ we 
essentially
have
$$
\|\widehat{f_\omega}\|_{L^2(B(x,R))}^2 \lesssim R^n X^* F_\omega(\frac{x}{R^2}).
$$
From this (and similar considerations for $g$) we see that \eqref{step-4}
is majorized by
 $$
 R^{-n} \int_{B(0,R^2)} (R^{-1/2} 
(\sum_{\omega \in \E_1} R^n X^* F_\omega(\frac{x}{R^2}))^{1/2} 
 R^{-1/2}
(\sum_{\omega \in \E_2} R^n X^* G_\omega(\frac{x}{R^2}))^{1/2}
	  )^q\ dx,
$$
where $G_\omega$ is defined in analogy to $F_\omega$.  We simplify this as
\begin{equation}\label{step-5}
R^{n+nq-q} \int_{B(0,1)} (X^* F(x) X^* G(x))^{q/2}\ dx,
\end{equation}
where
$$
F(\omega,i) = 
R^{n-1} \int_{\tilde T^\omega_i} |\widehat{\tilde f_\omega}(R^2 x)|^2\
dx, \quad
G(\omega,i) = 
R^{n-1} \int_{\tilde T^\omega_i} |\widehat{\tilde g_\omega}(R^2 x)|^2\
dx.$$
On the other hand, from the definition of the
hypothesis $K^*(\frac{p}{2} \times \frac{p}{2}
\to \frac{q}{2})$ we have
$$
\|  X^*_{1/R} F X^*_{1/R} G \|_{q/2}  \lesssim
R^{\frac{2n}{q^\prime}-2+\eps} 
\|F\|_{L^{p/2}_\omega L^1_i} \|G\|_{L^{p/2}_\omega L^1_i}.
$$
Comparing this with \eqref{step-3} and \eqref{step-5}, we see that
we will be done once we show that
$$ R^{\frac{2}{q}(n+nq-q)} R^{\frac{2n}{q^\prime}-2}
\|F\|_{L^{p/2}_\omega L^1_i} \|G\|_{L^{p/2}_\omega L^1_i}
\lesssim 
(R^{-2/\pp} \|f\|_p R^{-2/\pp} \|g\|_p)^2.$$
After some algebraic manipulation we see that it suffices to show that
$$ R^{2n-\frac{4}{p}+2} \|F\|_{L^{p/2}_\omega L^1_i} \lesssim \|f\|_{p}^2,$$
together with the completely analogous estimate for $g$, $G$.
From the definition of $f_\omega$ and the measure $d\omega$ we have
$$ \|f\|_p^2 \sim R^{\frac{n-1}{p/2}} 
\| ( \|f_\omega\|_p )^2 \|_{L^{p/2}_\omega},
$$
and so it suffices to show that
$$ R^{2n-\frac{4}{p}+2} \| F(\omega, \cdot) \|_{L^1_i} \lesssim R^{\frac{2(n-1)}{p}} 
\| f_\omega \|_{p}^2$$
uniformly in $\omega$.
From H\"older's inequality\footnote{This use of H\"older's inequality
indicates some room for improvement in this Lemma.  Indeed, one can
replace the $L^p$ norms on $f$, $g$ by the $B_p$ norms as used in
\cite{vargas:restrict}, Lemma 2.2.}, the hypothesis $p \geq 2$ and the support 
conditions on $f_\omega$ we have
$$ \| f_\omega \|_2 \lesssim R^{-(n+1)(\frac{1}{2}-\frac{1}{p})} \|f_\omega\|_p,$$
and so after some algebra we reduce ourselves to 
$$ \| F(\omega, \cdot) \|_{L^1_i} \lesssim R^{-n-1} \|f_\omega\|_2^2.$$
However, the left-hand side is majorized by
$$ R^{n-1} \int_{B(0,C)} |\widehat{\tilde f_\omega}(R^2 x)|^2\ dx
\lesssim R^{-n-1} \|\widehat{\tilde f_\omega}\|_2^2,$$
and the claim follows from Plancherel's theorem and the pointwise
comparability of $f_\omega$ and $\tilde f_\omega$.
\end{proof}

Applying Lemma \ref{b1} with $p=\frac{5}{2}$ and $q = \frac{5}{3}$ and using
\eqref{5-6}, we see that
$$ \R^*(2 \times 2 \to \frac{5}{3}, \alpha) \implies \R^*(\frac{5}{2}
\times \frac{5}{2} \to \frac{5}{3}, \frac{\alpha}{2}+\eps)$$
for all $\alpha,\eps > 0$.  On the other hand, from interpolating 
\eqref{trace} with \eqref{12-7} we obtain
$$
\R^*(\frac{30}{17} \times \frac{30}{17} \to \frac{5}{3}, \frac{1}{5}),
$$
so by another interpolation we obtain the implication
$$  \R^*(\frac{5}{2} \times \frac{5}{2} \to \frac{5}{3}, \beta)
\implies \R^*(2 \times 2 \to \frac{5}{3}, \frac{2}{5}\beta + \frac{3}{25})$$
for all $\beta > 0$.  Combining these two implications we see that
$$ \R^*(2 \times 2 \to \frac{5}{3}, \alpha) \implies
\R^*(2 \times 2 \to \frac{5}{3}, \frac{1}{5}\alpha + \frac{3}{25} + \eps).$$
The map $\alpha \to \frac{1}{5}\alpha + \frac{3}{25}$ is
a contraction with fixed point $\alpha = \frac{3}{20}$.  Since
the estimate $\R^*(2 \times 2 \to \frac{5}{3}, \alpha)$ holds for 
at least one value of $\alpha$, we thus see that
\begin{equation}\label{3-20}
 \R^*(2 \times 2 \to \frac{5}{3}, \frac{3}{20} + \eps)
 \end{equation}
for all $\eps > 0$.  Applying Lemma \ref{b1} one more time, we obtain
\begin{equation}\label{3-40}
 \R^*(\frac{5}{2} \times \frac{5}{2} \to \frac{5}{3}, \frac{3}{40} + \eps).
 \end{equation}

An inspection of the proof of Theorem \ref{imply} shows that the statement
of the theorem still holds when $\R^*(p \times p \to q)$ and $\R^*(p
\to 2q)$ are replaced by their local analogues
$\R^*(p \times p \to q,\alpha)$ and $\R^*(p \to 2q, \alpha/2)$.
Applying this to \eqref{3-40} we obtain
\begin{equation}\label{3-80}
 \R^*(\frac{5}{2} \to \frac{10}{3}, \frac{3}{80} + \eps).
\end{equation}

We now remove the $\alpha$ completely, borrowing the 
following argument of Bourgain
\cite{borg:kakeya,borg:stein} (for the concrete case
$n=3$, $p>20/7$, $q>10/3$, $\alpha>1/20$, $\ptil>7/3$, $\qtil>42/11$, see
 \cite{vargas:restrict}):

\begin{lemma}\label{alpha}\cite{borg:kakeya,borg:stein,vargas:restrict}
  If $p, q, \alpha$ are such that $\frac{n+1}{2} > \alpha q$,
then
$\R^*(p \to q, \alpha)$ implies $\R^*(\ptil \to \qtil)$ whenever
$$ \qtil > 2 + \frac{q}{\frac{n+1}{2} - \alpha q}, \quad
\frac{\qtil}{\ptil} <1 + \frac{q}{p} \frac{1}{\frac{n+1}{2}-\alpha q}.$$
\end{lemma}

Applying this to \eqref{3-80} we obtain the first conclusion of
Theorem \ref{main}.  Using Theorem \ref{imply} to return to the bilinear
setting, we thus obtain
$$ \R^*(p \times p \to q) \hbox{ for }
p > \frac{170}{77}, q > \frac{17}{9}.$$
Interpolating this with \eqref{12-7} and using Theorem \ref{imply},
one obtains the second conclusion of
Theorem \ref{main}.

\begin{figure}[htbp] \centering
\ \psfig{figure=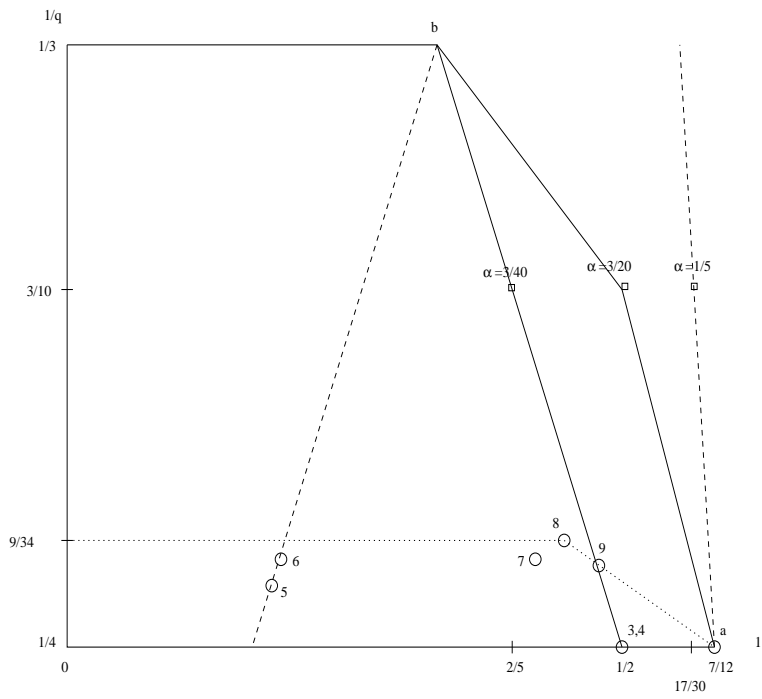,height=3.6in,width=4in}
\caption{Estimates of the form $R^*(p \times p \to q,\alpha)$ and 
$R^*(p \to 2q,\alpha/2)$ for $n=3$.  }
        \label{fig:rest}
	\end{figure}

We summarize the various estimates used in Figure \ref{fig:rest},
which is an expanded version of Figure \ref{fig:bilinear}.
For comparison, the previously known results are also displayed.  The
dotted line thus represents the best global restriction theorems
(both linear and bilinear) known to date.  (It is possible to
improve on these results slightly; see \cite{mtvv:cone}).

By interpolating between the main result and \eqref{12-7} we also
obtain some progress on Klainerman's conjecture for the sphere in $\R^3$:

\begin{corollary}\label{beta} If $n=3$, then $R^*(2 \times 2 \to p)$ whenever
$p > 2 - \frac{5}{69}$.
\end{corollary}

These techniques are certainly not best possible.  For instance,
one can use the techniques in \cite{borg:cone} to obtain better
versions of Corollary \ref{beta}. See \cite{mtvv:cone}.

The sharp restriction theorem $R^*_s(q)$ is scale-invariant under parabolic
scaling.  Thus, the compact support condition on $\Phi$ can be removed.
In particular, one has a sharp restriction theorem for the entire paraboloid
$\{(\underline{x}, \frac{1}{2}|\underline{x}|^2): 
\underline{x} \in \R^{n-1}\}$ for $q > 4 - \frac{5}{27}$.

One can extend the above results to Bochner-Riesz
multipliers, so that the Bochner-Riesz conjecture
holds for $n=3$ and $\max(p,\pp) \geq \frac{34}{9}$.  We sketch
the argument very briefly as follows.  By the
usual techniques of Carleson-Sj\"olin reduction and factorization
theory (see \cite{borg:stein}) it suffices to show that
$$ \| Tf\|_{L^p(\R^n)}
\lesssim \lambda^{-n/p} \|f\|_{L^\infty(Q)}$$
for all $p > 34/9$, $\lambda \gg 1$ and $f \in L^\infty(Q)$, where 
$$ Tf(x) = \int_Q e^{2\pi i \lambda|x-y|} a(x,y) f(y)\ dy,$$
$Q$ is thought of as imbedded in $\R^n$, 
and $a$ is a bump 
function on $\R^n \times Q$ which is supported away from the diagonal
$x=y$.  By the analogue of Lemma \ref{alpha} for Bochner-Riesz multipliers
(see \cite{borg:stein}) it suffices to show that
$$ \| Tf \|_{10/3}
\lesssim \lambda^{3/80+\eps} \lambda^{-9/10} \|f\|_\infty$$
for all $\eps > 0$.  By a modification of Theorem \ref{imply} it suffices
to show that
$$ \| Tf Tg\|_{5/3} \lesssim \lambda^{3/40+\eps} \lambda^{-9/10} \|f\|_\infty
 \lambda^{-9/10} \|g\|_\infty$$
for all $f$, $g$ with $O(1)$ separated supports, together with variants
of this estimate in which the phase function $|x-y|$ is replaced by
a parabolically scaled (but essentially equivalent) version.  
However, from the analogue
of Lemma \ref{b1} for Bochner-Riesz operators (which is proven
by a bilinear modification of the arguments in \cite{borg:stein})
this will follow from the restriction estimate \eqref{3-20}
and the analogue of Theore \ref{b-wolff} for the Nikodym maximal operator
(see e.g. \cite{wolff:kakeya}), which is proven similarly.  The 
required Nikodym estimate also follows formally from the original
formulation of Theorem \ref{b-wolff}: see the argument in \cite{tao:boch-rest}.

In higher dimensions $n>3$ Theorem \ref{modest} becomes too weak to
be of much use,
and we can only achieve a minor improvement
on known results.  By interpolating between \eqref{trace} and
 the bilinear form
$ R^*(2 \times 2 \to \frac{n+1}{n-1})$ of the Tomas-Stein theorem,
we obtain
$$ R^*(2 \times 2 \to \frac{n+2}{n}, \frac{1}{n+2}).$$
Applying this and Theorem \ref{b-wolff} to Lemma \ref{b1} we obtain
$$ R^*(\frac{2(n+2)}{n+1} \times \frac{2(n+2)}{n+1} \to
\frac{n+2}{n}, \frac{1}{2(n+1)} + \eps).$$
Applying Theorem \ref{imply} this becomes
$$ R^*(\frac{2(n+2)}{n+1} \to \frac{2(n+2)}{n}, \frac{1}{4(n+1)} + \eps).$$
Applying Lemma \ref{alpha} this becomes
$$ R^*(p \to q) \hbox{ for } p > \frac{2n^2 + 6n + 6}{n^2 + 3n + 1},
\quad q > \frac{2n^2 + 6n + 6}{n^2 + n - 1}.$$
This is only a slight improvement on the result in Wolff \cite{wolff:kakeya},
which showed $R^*(q \to q)$ for the same range of $q$.  For $n>3$ the
results obtained by interpolating these estimates with Theorem
\ref{modest} are inferior to the Tomas-Stein theorem.

\section{Further remarks}

In the previous sections we obtained a non-trivial
sharp restriction theorem $R^*_s(4-\eps)$ from an ordinary restriction
theorem $R^*(p \to q)$ (in this case $p > \frac{170}{77}, q >  
\frac{34}{9}$)
and the bilinear estimate \eqref{12-7}.

The original formulation of \eqref{12-7} in \cite{vargas:restrict,vargas:2} 
was stated in
terms of $X_r$ spaces.  In this section we show how one can use these
estimates instead of the bilinear estimate to obtain non-trivial sharp
restriction theorems.  Despite the fact that these estimates can
be extended (for characteristic functions) from $r > \frac{12}{7}$
to $r \geq 4(\sqrt{2}-1)$, the methods we will use do not appear
to be as efficient as the bilinear techniques.
However, they seem to be more robust and applicable to a wider
range of situations.

\begin{proposition}\label{x-imply}  Let $n=3$.  
Suppose that $R^*(p \to q)$ holds for some
$2 < q < 4$, and suppose that the quantity
$r = \frac{4\pp}{q}$ satisfies $r> 4(\sqrt{2}-1)$.
Then we have $R^*_s(w)$ for all $w > \frac{4+q}{2}$.
\end{proposition}

Note that
the points $(\frac{1}{p},\frac{1}{q}), (1-\frac{2}{w},\frac{1}{w}), 
(\frac{1}{r},\frac{1}{4})$ are collinear when $w = \frac{4+q}{2}$.

\begin{proof}
It suffices to show the restricted weak-type estimate
\begin{equation}\label{w-est}
|\{|\Re^* \chi_\Omega| \gtrsim \lambda\}| \lesssim \frac{2^{-2(w-2)j_0}}{\lambda^w}
\end{equation}
for all $\lambda > 0$ and $\Omega \subseteq Q$, where $j_0$ is
the integer such that $|\Omega| \sim 2^{-2j_0}$.  
We may assume that $2^{-C j_0} \lesssim \lambda
\lesssim 1$ for some constant $C$, since \eqref{w-est} is trivial
(by e.g. the Tomas-Stein restriction theorem) otherwise.

The idea of the proof will be to decompose $\Omega$ into a sparse set
and a collection of sets concentrated on caps.  
On the sparse set the Tomas-Stein estimate $R^*_s(4)$ can be improved
using the $X_r$ estimates of \cite{vargas:restrict,vargas:2}.
The sets on caps
can be rescaled parabolically to become sets of measure comparable to 1,
in which case the estimate $R^*(p \to q)$ is equivalent to the
scale-invariant estimate $R^*_s(q)$.  Combining the two estimates one expects
to obtain $R^*_s(4-\eps)$ for some $\eps > 0$.

We now turn to the details.  
For each $j$ let $0 < \alpha_j \leq 1$ be a quantity to be chosen later.
By the usual Calder\'on-Zygmund stopping time arguments, we may 
partition $\Omega$ as
$$\Omega= \Omega_g \cup \bigcup_{(j,k) \in T} (\Omega \cap \tau^j_k),$$
where the ``good'' set $\Omega_g$ satisfies
\begin{equation}\label{good-def}
 |\Omega_g \cap \tau^j_k| \leq \alpha_j |\tau^j_k|
\end{equation}
for all $j,k$, and $\{ \tau_{j,k}: (j,k) \in T \}$ is a collection of 
disjoint dyadic cubes such that
\begin{equation}\label{bad}
\alpha_j |\tau^j_k| < |\Omega \cap \tau^j_k| \leq 4 \alpha_{j-1} |\tau^j_k|
\end{equation}
for all $(j,k) \in T$.

We decompose $\Re^* \chi_\Omega$ as
$$ \Re^* \chi_\Omega = \Re^* \chi_{\Omega_g} + \sum_j \sum_{k: (j,k) \in T}
\Re^* \chi_{\Omega \cap \tau^j_k}.$$

To control the contribution of $\Re^* \chi_{\Omega_g}$ we use
the results of \cite{vargas:restrict,vargas:2} and the hypothesis
$r > 4(\sqrt{2}-1)$ to obtain 
\begin{equation}\label{xr}
\|\,\Re^* \chi_{\Omega_g}\,\|_4\lesssim \|\chi_{\Omega_g}\|_{X_r}
= \bigl(\sum_j\sum_k 2^{-4j}
\bigl(\,\frac{|\Om_g \cap\tau^j_k|}{|\tau^j_k|}\bigr)^{4/r}\bigr)^{1/4}.
\end{equation}

In order for $R^*(p \to q)$ to hold we must have $p^\prime \leq \frac{1}{2}q$,
so that $4/r \geq 2 > 1$.  Applying \eqref{xr-big} of Lemma \ref{xr-est} 
with $p=4/r$ and $\alpha = \alpha_j$ we obtain (after some algebra)
 $$
 \sum_k 2^{-4j} \bigl(\,\frac{|\Om_g\cap\tau^j_k|}{|\tau^j_k|}\bigr)^{4/r}
\lesssim \min(
2^{(\frac{8}{r} - 4)(j-j_0)} 2^{-4j_0},
\alpha_j^{\frac{4}{r}-1} 2^{2(j_0-j)} 2^{-4j_0} 
)
$$
for each $j$; informally, this shows that the most significant scales
occur when $j$ is near $j_0$.  Inserting these estimates into \eqref{xr} we 
obtain

\begin{equation}\label{om-g}
\|\,\Re^* \chi_{\Omega_g}\,\|_4^4\lesssim 
\sum_j 
\min(
2^{(\frac{8}{r} - 4)(j-j_0)},
\alpha_j^{\frac{4}{r}-1}
2^{2(j_0-j)}) 2^{-4j_0}.
\end{equation}
Let $m$ be a positive integer to be chosen later.  We now choose
$\alpha_j$ so that
\begin{equation}\label{alpha-def}
\alpha^{\frac{4}{r}-1} 2^{2(j_0-j)}
= 2^{-m}
\end{equation}
when $-\frac{m}{\frac{8}{r}-4} < j-j_0 < \frac{m}{2}$, and
$\alpha_j = 1$ otherwise.  The estimate \eqref{om-g} then becomes
$$ \|\,\Re^* \chi_{\Omega_g}\,\|_4^4\lesssim m 2^{-m} 2^{-4j_0}.$$
From Tchebyshev's inequality we thus obtain
\begin{equation}\label{good-t}
 | \{\Re^* \chi_{\Omega_g} \gtrsim \lambda \} | \lesssim 
m 2^{-m} 2^{-4j_0} \lambda^{-4}.
\end{equation}

It remains to estimate the quantity
\begin{equation}\label{bad-part}
\sum_{k: (j,k) \in T} \Re^* \chi_{\Omega \cap \tau^j_k}
\end{equation}
for each $j$.  Note that this quantity vanishes when $\alpha_j=1$
by \eqref{bad}, so we may assume that $-\frac{m}{\frac{8}{r}-4} < j-j_0 < 
\frac{m}{2}$.

We will estimate \eqref{bad-part} in $L^q$ norm.  From Lemma \ref{young}
in the Appendix we have
\begin{equation}\label{quasi}
 \| \sum_{k: (j,k) \in T} \Re^* \chi_{\Omega \cap \tau^j_k} \|_q
\lesssim (\sum_{k: (j,k) \in T} \|\Re^* \chi_{\Omega \cap \tau^j_k} \|_q^\qp)^{1/\qp}.
\end{equation}
On the other hand, by parabolically rescaling the hypothesis $R^*(p \to q)$ 
as in Proposition \ref{rescale} we obtain
$$
\|\Re^* \chi_{\Om \cap \tau^j_k} \|_{L^q}\lesssim
|\tau^j_k|^{1 - \frac{2}{q}-\frac{1}{p}} |\Omega \cap \tau^j_k|^{1/p}.
$$
Combining this with \eqref{quasi} and \eqref{bad} we obtain
$$ \| \sum_{k: (j,k) \in T} \Re^* \chi_{\Omega \cap \tau^j_k} \|_q
\lesssim \bigl(\sum_{k: (j,k) \in T} 
(|\tau^j_k|^{1 - \frac{2}{q}} (4\alpha_{j-1})^{1/p})^\qp\bigr)^{1/\qp}.$$
On the other hand, from \eqref{bad} we have
$$ \# \{k: (j,k) \in T\} \leq 2^{2(j-j_0)} \alpha_j^{-1}.$$
Using this and the fact that $4 \alpha_{j-1} \sim \alpha_j$ we obtain
\begin{equation} \label{final-j}
 \| \sum_{k: (j,k) \in T} \Re^* \chi_{\Omega \cap \tau^j_k} \|_q
\lesssim |\tau^j_k|^{1 - \frac{2}{q}} 2^{\frac{2}{\qp}(j-j_0)} 
\alpha_j^{\frac{1}{p} - \frac{1}{\qp}}.
\end{equation}
Combining this with \eqref{alpha-def} and the definition of $r$
one eventually obtains
$$ \| \sum_{k: (j,k) \in T} \Re^* \chi_{\Omega \cap \tau^j_k} \|_q
\lesssim 2^{-2(1-\frac{2}{q})j_0} 2^{\frac{m}{q}}$$
so by the triangle inequality and Tchebyshev's inequality we have
$$ \{ |\sum_j \sum_{k: (j,k) \in T} \Re^* \chi_{\Omega \cap \tau^j_k}|
\gtrsim \lambda \} \lesssim 2^{-2(q-2)j_0} m^q 2^{m} \lambda^{-q}.$$
Combining this with \eqref{good-t} we obtain
$$ \{ |\Re^* \chi_\Omega| \gtrsim \lambda \} \lesssim
m 2^{-m} 2^{-4j_0} \lambda^{-4} + 2^{-2(q-2)j_0} m^q 2^{m} 
\lambda^{-q}.$$
The claim \eqref{w-est} then follows by choosing 
$2^m = (2^{2j_0} \lambda)^{-\frac{4-q}{2}}$.
\end{proof}

In practice Proposition \ref{x-imply} is inferior to the implications obtained 
by Theorem \ref{imply} and interpolation with \eqref{12-7}.  For instance,
if we insert the first conclusion of Theorem \ref{main} into Proposition
\ref{x-imply} one obtains $R^*_s(w)$ for $w > 4 - \frac{1}{9}$, which 
is inferior to the second conclusion of Theorem \ref{main}.

\section{Appendix: Some elementary harmonic analysis}

In this section we state some elementary results which were used
repeatedly in the paper.  

We begin with a well-known quasi-orthogonality
property of functions with disjoint frequency support.  Define 
a \emph{rectangle} to be the product of $n$ (possibly half-infinite
or infinite) intervals in $\R^n$.

\begin{lemma}\label{quasi-lemma} Let $R_k$ be a collection of rectangles
in frequency space such that the dilates $2R_k$ are almost disjoint,
and suppose that $f_k$ are a collection of functions whose Fourier
transforms are supported on $R_k$.  Then for all $1 \leq p \leq \infty$
we have
$$ \| \sum_k f_k \|_p \lesssim (\sum_k \|f_k\|_p^{p^*})^{1/p^*},$$
where $p^* = \min(p, p^\prime)$.
\end{lemma}

\begin{proof}  Let $P_k$ be a smooth Fourier multiplier adapted to $2R_k$
which equals $1$ on $R_k$.  We claim that
$$ \| \sum_k P_k F_k \|_p \lesssim (\sum_k \|F_k\|_p^{p^*})^{1/p^*}$$
for arbitrary functions $F_k$; the lemma then follows by setting
$F_k = P_k F_k = f_k$.

By interpolation it suffices to prove this estimate for $p = 1$, $p=2$, 
and $p=\infty$.  When $p=2$ the estimate is immediate from Plancherel's
theorem.  When $p=1$ or $p=\infty$ the lemma follows from the triangle
inequality and the estimates
$$ \| P_k F_k \|_1 \lesssim \|F_k\|_1, \quad 
\| P_k F_k \|_\infty \lesssim \|F_k\|_\infty,$$
which follow from Young's inequality and standard estimates on the kernel
of $P_k$.
\end{proof}

The next lemma allows us to crudely estimate various $X_r$-type quantities.

\begin{lemma}\label{xr-est} 
Let $\Omega \subseteq Q$ be a set such that $|\Omega| \lesssim
2^{-(n-1)j_0}$ for some $j_0 \geq 0$, and let $\tau^j_k$ be defined as in
Section \ref{iff}.  Let $0 \leq \alpha \leq 1$ be such that
$|\Omega \cap \tau^j_k| \leq \alpha |\tau^j_k|$ for all $k$.
If $p \geq 1$, then we have the estimate
\begin{equation}\label{xr-big}
 \sum_k |\Omega \cap \tau^j_k|^p \lesssim 2^{-(n-1) j_0} 
 \min(\alpha 2^{-(n-1)j}, 2^{-(n-1)j_0})^{p-1}.
\end{equation}
If $p \leq 1$, then we have the estimate
\begin{equation}\label{xr-small}
 \sum_k |\Omega \cap \tau^j_k|^p \lesssim 2^{-(n-1)p j_0}
 2^{(n-1)(1-p) j}.
\end{equation}
\end{lemma}

\begin{proof}
By the log-convexity of $L^p$ norms, $0 \leq p \leq \infty$, it suffices to 
prove the three bounds
\begin{align}
\sup_k |\Omega \cap \tau^j_k| &\lesssim 
 \min(\alpha 2^{-(n-1)j}, 2^{-(n-1)j_0})
\label{x1}\\
\sum_k |\Omega \cap \tau^j_k| &\lesssim 2^{-(n-1) j_0}\label{x2}\\
\sum_k 1 &\lesssim 2^{(n-1)j}.\label{x3}
\end{align}
But these bounds follow trivially from the estimates
$$ |\Omega \cap \tau^j_k| \leq \min(|\Omega|,\alpha |\tau^j_k|) \lesssim \min(2^{-(n-1)j_0}, \alpha 2^{-(n-1)j}),$$
$$ \sum_k |\Omega \cap \tau^j_k| = |\Omega| \lesssim 2^{-(n-1)j_0},$$
and the cardinality of the $\tau^j_k$.
\end{proof}

We remark that without further information on $\Omega$ these bounds
are best possible.  In most cases we will set $\alpha=1$.

Finally, we present a very easy inequality.

\begin{lemma}\label{young}  If $1 \leq p \leq \infty$, then
$$ (\sum_k |a_k|^p)^{1/p} \leq \sum_k |a_k|$$
for all sequences of numbers $a_k$.  Also, if $0 < q \leq 1$, then
$$ \|\sum_k f_k\|_q \leq (\sum_k \|f_k\|_q^q)^{1/q}.$$
\end{lemma}

\begin{proof}
The first estimate is trivial for $p=1$ and $p=\infty$, and the
general case follows by convexity.  The second estimate follows
by applying the first with $a_k = |f_k(x)|^q$, $p = 1/q$ and integrating.
\end{proof}


\begin{thebibliography}{10}

\bibitem{beals:xsb}
M. Beals, \emph{Self-Spreading and strength of Singularities for solutions to
semilinear wave equations}, Annals of Math \textbf{118} (1983), 187-214.
  
\bibitem{borg:kakeya}
J. Bourgain, \emph{Besicovitch-type maximal operators and
applications to Fourier analysis}, Geom. and Funct. Anal. \textbf{22}
(1991), 147--187.

\bibitem{borg:16-9}
J. Bourgain, \emph{On the restriction and multiplier problem
in $\R^3$}, Lecture notes in Mathematics, no. 1469.  Springer Verlag, 1991.

\bibitem{borg:schrodinger}
J. Bourgain, \emph{A remark on Schrodinger operators},
Israel J. Math. 77 (1992), 1--16.

\bibitem{borg:cone}
J. Bourgain, \emph{Estimates for cone multipliers},
Operator Theory: Advances and Applications, \textbf{77} (1995), 41--60.

\bibitem{borg:stein}
J. Bourgain, \emph{Some new estimates on oscillatory integrals},
Essays in Fourier Analysis in honor of E.~M. Stein, Princeton University
Press (1995), 83--112.

\bibitem{carl:disc}
L. Carleson and P. Sj\"olin, \emph{Oscillatory integrals and a
multiplier problem for the disc}, Studia Math. \textbf{44} (1972): 287--299.

\bibitem{cordoba:sieve}
A. C\'ordoba, \emph{The Kakeya maximal function and the spherical summation
multipliers}, Amer. J. Math. 99 (1977), 1--22.

\bibitem{drury:xray}
S. Drury, \emph{$L^p$ estimates for the x-ray transform}, Ill. J. Math. 
\textbf{27} (1983), 125--129.

\bibitem{feff:thesis}
C. Fefferman, \emph{Inequalities for strongly
singular convolution operators}, Acta Math. \textbf{ 124}
(1970), 9--36.

\bibitem{feff:ball}
C. Fefferman, \emph{The multiplier problem for the ball}, Ann. of Math. \textbf{
94} (1971), 330--336.

\bibitem{horm:fio}
L. H\"ormander, \emph{Fourier Integral Operators},
Acta Math. \textbf{127} (1971): 79--183.

\bibitem{katz:kakeya}
N. Katz, preprint.

\bibitem{kl-mac:cpde-null}
S. Klainerman, M. Machedon, \emph{Space-time estimates for null forms
and the local existence theorem}, Comm. Pure Appl. Math. 46 (1993), no. 9,
1221--1268. 

\bibitem{kl-mac:null}
S. Klainerman, M. Machedon, \emph{Remark on Strichartz-type
inequalities.} With appendices by Jean Bourgain and Daniel Tataru. 
Internat. Math.  Res. Notices \textbf{5} (1996), 201--220.

\bibitem{kl-mac:duke-null}
S. Klainerman, M. Machedon, \emph{On the regularity properties
of a model problem related to wave maps}, Duke Math. J. 87 (1997),
553--589.

\bibitem{vargas:restrict}
A. Moyua, A. Vargas, L. Vega, \emph{Schr\"odinger Maximal Function
and Restriction Properties of the Fourier transform}, International Math.
Research Notices \textbf{16} (1996).

\bibitem{vargas:2}
A. Moyua, A. Vargas, L. Vega, \emph{Restriction theorems and Maximal
operators related to oscillatory integrals in $\R^3$}, to appear, Duke
Math. J.

\bibitem{schlag:bourgain}
W. Schlag, \emph{A generalization of Bourgain's circular maximal
theorem}, J. Amer. Math. Soc. \textbf{10} (1997), 103-122.

\bibitem{schlag:geom-borg}
W. Schlag, \emph{A geometric proof of the circular maximal theorem},
to appear, Duke Math. J.

\bibitem{schlag:kakeya}
W. Schlag, \emph{A geometric inequality with applications to the
Kakeya problem in three dimensions}, 
to appear, Geometric and Functional Analysis.

\bibitem{sogge:kakeya}
C.~D. Sogge, \emph{Concerning Nikodym-type sets in 3-dimensional
curved space}, preprint.

\bibitem{stein:large}
E.~M. Stein, \emph{Harmonic Analysis}, Princeton University Press, 1993.

\bibitem{tao:boch-rest}
T. Tao, \emph{The Bochner-Riesz conjecture implies the Restriction
conjecture}, to appear, Duke Math J.

\bibitem{mtvv:cone}
T. Tao, A. Vargas, \emph{A bilinear approach to cone multipliers and
related operators}, in preparation.

\bibitem{tomas:restrict}
P. Tomas, \emph{A restriction theorem for the Fourier transform}, Bull. Amer.
Math. Soc. \textbf{81} (1975), 477--478.

\bibitem{wolff:kakeya}
T.~H. Wolff, \emph{An improved bound for Kakeya type maximal functions},
Revista Mat. Iberoamericana. \textbf{11} (1995). 651--674.

\bibitem{wolff:xray}
T.~H. Wolff, \emph{A mixed norm estimate for the x-ray transform}, to 
appear in Revista Mat. Iberoamericana.



\end{thebibliography}
\end{document}